# Multivariate COGARCH(1, 1) processes

ROBERT STELZER

*TUM Institute for Advanced Study & Zentrum Mathematik, Technische Universität München, Boltzmannstraße 3, D-85747 Garching, Germany. E-mail: rstelzer@ma.tum.de*

Multivariate COGARCH(1, 1) processes are introduced as a continuous-time models for multidimensional heteroskedastic observations. Our model is driven by a single multivariate Lévy process and the latent time-varying covariance matrix is directly specified as a stochastic process in the positive semidefinite matrices.

After defining the COGARCH(1, 1) process, we analyze its probabilistic properties. We show a sufficient condition for the existence of a stationary distribution for the stochastic covariance matrix process and present criteria ensuring the finiteness of moments. Under certain natural assumptions on the moments of the driving Lévy process, explicit expressions for the first and second-order moments and (asymptotic) second-order stationarity of the covariance matrix process are obtained. Furthermore, we study the stationarity and second-order structure of the increments of the multivariate COGARCH(1, 1) process and their "squares".

## 1. Introduction

In this paper, a multivariate extension of the continuous-time generalized autoregressive conditional heteroskedasticity (COGARCH, for short) process of order (1, 1) introduced in [26] is defined and studied in detail. The one-dimensional COGARCH(1, 1) process (see also [10, 22, 27]) is given as the solution of

$$\mathrm{d}G_t = \sqrt{v_{t-}}\,\mathrm{d}L_t, \tag{1.1}$$

$$\mathrm{d}v_t = -\beta(v_{t-} - c)\,\mathrm{d}t + \alpha v_{t-}\,\mathrm{d}[L, L]_t^{\mathfrak{d}}, \tag{1.2}$$

using the discontinuous part $[L, L]^{\mathfrak{d}}$ of the quadratic variation of a univariate Lévy process $L$, parameters $\alpha, \beta, c > 0$ and initial values $G_0 = 0, v_0 \geq 0$. The process $G$ is referred to as the COGARCH(1, 1) *process* and the variance process $v$ as its *volatility process*, where the name "volatility process" derives from the term typically used in economics.

Heteroskedastic data are often modelled with (normal) variance mixture models. In such a model, one has $X_n = \sqrt{v_n}\varepsilon_n$ for $n \in \mathbb{N}$, where $\varepsilon$ is an i.i.d. sequence and $v$ a sequence of positive random variables modelling the current variance of the observations $X$. Typically, one has, moreover, that $\varepsilon_n$ and $v_n$ are, for fixed $n \in \mathbb{N}$, independent. Obviously, equations (1.1) and (1.2) constitute a continuous-time counterpart of a variance mixture model with some special process for the variance and driven by a single Lévy process. Loosely speaking, the increments $\mathrm{d}L_t$ are "mixed" with the variance $v_{t-}$ and the two are independent for fixed $t$.

In a multivariate setup, the positive variance $v$ needs to be replaced by a covariance matrix process $V$. Thus, the volatility process has to be a stochastic process in the positive semidefinite matrices. This requirement leads to challenging questions in modelling and interesting mathematical issues since, in particular, only very few continuous-time stochastic processes in the positive semidefinite matrices have been thus far studied (mainly various "Wishart" processes, see [11, 12, 14, 19, 20]; or, recently, Ornstein–Uhlenbeck-type processes, see [5, 34, 35]). Appropriate multivariate models for heteroskedastic data are, however, clearly needed because in many areas of application, one has to model and understand the joint behavior of several time series exhibiting nontrivial interdependencies. Moreover, for various reasons (for example, unequally spaced observations, inference at several frequencies or amenability to continuous-time financial theory), it is often desirable to use continuous-time models instead of related discrete-time models like GARCH models, for instance.

After briefly stating some preliminaries regarding notation and Lévy processes in Section 2, we introduce our multivariate COGARCH(1, 1) processes (MUCOGARCH(1, 1), for short) in Section 3 and establish well-definedness. Thereafter, we analyze its volatility process in Section 4. In the first part of that section, we present a univariate COGARCH(1, 1) process that bounds the volatility process in a norm intrinsically related to the autoregressive parameter and use this bound to give sufficient conditions for the finiteness of moments. This is followed by a demonstration that the volatility process alone and the MUCOGARCH(1, 1) process together with its volatility are strong Markov processes. Moreover, we establish conditions for the existence of a stationary distribution of the volatility in Section 4.2. In the last part of Section 4, we calculate the second-order structure of the volatility process explicitly under certain assumptions on the moments of the driving Lévy process and establish (asymptotic) second-order stationarity.

In Section 5, we focus on the increments of the MUCOGARCH(1, 1) process itself, showing that it has stationary increments provided the volatility is stationary. Thereafter, we calculate the second-order moment structure of the increments (that is, the returns in a financial context) observed on a regularly spaced discrete grid and their "squares" (that is, the increments times their transposes). Here, we obtain, in particular, that the increments have zero autocorrelation, but their "squares" have exponentially decaying autocorrelation. Moreover, the explicit expressions for the moments obtained make the processes amenable to statistical estimation.

Finally, we present all proofs, together with auxiliary technical results, in Section 6.

# 2. Preliminaries

## 2.1. Notation

We denote the set of real $m \times n$ matrices by $M_{m,n}(\mathbb{R})$. If $m = n$, we simply write $M_n(\mathbb{R})$ and denote the group of invertible $n \times n$ matrices by $GL_n(\mathbb{R})$, the linear subspace of symmetric matrices by $\mathbb{S}_n$, the (closed) positive semidefinite cone by $\mathbb{S}_n^+$ and the open positive definite cone by $\mathbb{S}_n^{++}$. $I_n$ stands for the $n \times n$ identity matrix. The natural ordering on the symmetric $n \times n$ matrices shall be denoted by $\leq$, that is, for $A, B \in M_n(\mathbb{R})$, we have that $A \leq B$ if and only if $B - A \in \mathbb{S}_n^+$ (likewise, $A < B \Leftrightarrow B - A \in \mathbb{S}_d^{++}$). The tensor (Kronecker) product of two matrices $A, B$ is written as $A \otimes B$. vec denotes the well-known vectorization operator that maps the set of $n \times n$ matrices to $\mathbb{R}^{n^2}$ by stacking the columns of the matrices below one another. For more information regarding the tensor product and vec operator, we refer to [24], Chapter 4. The spectrum of a matrix is denoted by $\sigma(\cdot)$ and the spectral radius by $\rho(\cdot)$. Finally, $A^*$ denotes the transpose (adjoint) of a matrix $A \in M_{m,n}(\mathbb{R})$.

Norms of vectors or matrices are denoted by $\|\cdot\|$. If the norm is not further specified, then it is irrelevant which particular norm is used.

Throughout, we assume that all random variables and processes are defined on a given filtered probability space $(\Omega, \mathcal{F}, P, (\mathcal{F}_t)_{t\in\mathcal{T}})$, with $\mathcal{T} = \mathbb{N}$ in the discrete-time case and $\mathcal{T} = \mathbb{R}^+$ in the continuous-time one. Moreover, in the continuous-time setting, we assume the usual conditions (complete, right-continuous filtration) to be satisfied.

Furthermore, we employ an intuitive notation with respect to (stochastic) integration with matrix-valued integrators, referring to any of the standard texts (for example, [36]) for a comprehensive treatment of the theory of stochastic integration. Let $(A_t)_{t\in\mathbb{R}^+}$ in $M_{m,n}(\mathbb{R})$ and $(B_t)_{t\in\mathbb{R}^+}$ in $M_{r,s}(\mathbb{R})$ be càdlàg and adapted processes and $(L_t)_{t\in\mathbb{R}^+}$ in $M_{n,r}(\mathbb{R})$ be a semimartingale. We then denote by $\int_0^t A_{s-}\, \mathrm{d}L_s B_{s-}$ the matrix $C_t$ in $M_{m,s}(\mathbb{R})$ which has $ij$th element $C_{ij,t} = \sum_{k=1}^n \sum_{l=1}^r \int_0^t A_{ik,s-} B_{lj,s-}\, \mathrm{d}L_{kl,s}$. Equivalently, such an integral can be understood in the sense of [32, 33] by identifying it with the integral $\int_0^t \mathbf{A}_{s-}\, \mathrm{d}L_s$, with $\mathbf{A}_t$ being, for each fixed $t$, the linear operator $M_{n,r}(\mathbb{R}) \to M_{m,s}(\mathbb{R}), X \mapsto A_t X B_t$. If $(X_t)_{t\in\mathbb{R}^+}$ is a semimartingale in $\mathbb{R}^m$ and $(Y_t)_{t\in\mathbb{R}^+}$ one in $\mathbb{R}^n$, then the quadratic variation $([X,Y]_t)_{t\in\mathbb{R}^+}$ is defined as the finite variation process in $M_{m,n}(\mathbb{R})$ with components $[X,Y]_{ij,t} = [X_i, Y_j]_t$ for $t \in \mathbb{R}^+$ and $i = 1, \ldots, m$, $j = 1, \ldots, n$.

## 2.2. Lévy processes

Later, we shall use Lévy processes (see [1, 36, 39], for instance) both in $\mathbb{R}^d$ and in the symmetric matrices $\mathbb{S}_d$. Thus, we briefly recall the relevant basic notions for them now.

We consider a Lévy process $L = (L_t)_{t\in\mathbb{R}^+}$ (where $L_0 = 0$ a.s.) in $\mathbb{R}^d$ which is determined by its characteristic function in the Lévy–Khintchine form $E[\mathrm{e}^{\mathrm{i}\langle u, L_t\rangle}] = \exp\{t\psi_L(u)\}$ for

$t \in \mathbb{R}^+$ with

$$\psi_L(u) = \mathrm{i}\langle\gamma_L, u\rangle - \frac{1}{2}\langle u, \tau_L u\rangle + \int_{\mathbb{R}^d} (\mathrm{e}^{\mathrm{i}\langle u,x\rangle} - 1 - \mathrm{i}\langle u, x\rangle I_{[0,1]}(\|x\|))\nu_L(\mathrm{d}x), \qquad u \in \mathbb{R}^d,$$

where $\gamma_L \in \mathbb{R}^d$, $\tau_L \in \mathbb{S}_d^+$ and the Lévy measure $\nu_L$ is a measure on $\mathbb{R}^d$ satisfying $\nu_L(\{0\}) = 0$ and $\int_{\mathbb{R}^d}(\|x\|^2 \wedge 1)\nu_L(\mathrm{d}x) < \infty$. Moreover, $\langle\cdot,\cdot\rangle$ denotes the usual Euclidean scalar product on $\mathbb{R}^d$.

We always assume $L$ to be cadlag and denote its jump measure by $\mu_L$, that is, $\mu_L$ is the Poisson random measure on $\mathbb{R}^+ \times \mathbb{R}^d \setminus \{0\}$ given by $\mu_L(B) = \sharp\{s \geq 0 : (s, L_s - L_{s-}) \in B\}$ for any measurable set $B \subset \mathbb{R}^+ \times \mathbb{R}^d \setminus \{0\}$. Likewise, $\tilde{\mu}_L(\mathrm{d}s, \mathrm{d}x) = \mu_L(\mathrm{d}s, \mathrm{d}x) - \mathrm{d}s\, \nu_L(\mathrm{d}x)$ denotes the compensated jump measure.

If $\int_{\|x\|>1} \|x\|^2 \nu_L(\mathrm{d}x) < \infty$, then $L$ has finite mean and covariance matrix given by

$$E(L_1) = \gamma_L + \int_{\|x\|>1} x\nu(\mathrm{d}x), \qquad \mathrm{var}(L_1) = \tau_L + \int_{\mathbb{R}^d} xx^* \nu_L(\mathrm{d}x). \tag{2.1}$$

Provided $\tau_L = 0$ and $\int_{\|x\|\leq 1} \|x\|\nu_L(\mathrm{d}x) < \infty$, the Lévy process $L$ has paths of finite variation and

$$\psi_L(u) = \mathrm{i}\langle\tilde{\gamma}_L, u\rangle + \int_{\mathbb{R}^d} (\mathrm{e}^{\mathrm{i}\langle u,x\rangle} - 1)\nu_L(\mathrm{d}x) \qquad \text{with } \tilde{\gamma}_L = \gamma_L - \int_{\|x\|\leq 1} x\nu_L(\mathrm{d}x).$$

Regarding matrix-valued Lévy processes, we will only encounter matrix subordinators (see [3]), that is, Lévy processes with paths in $\mathbb{S}_d^+$. Since matrix subordinators are of finite variation and $\mathrm{tr}(X^*Y)$ (with $X, Y \in \mathbb{S}_d$ and tr denoting the usual trace functional) defines a scalar product on $\mathbb{S}_d$ linked to the Euclidean scalar product on $\mathbb{R}^{d^2}$ via $\mathrm{tr}(X^*Y) = \mathrm{vec}(X)^* \mathrm{vec}(Y) = \langle \mathrm{vec}(Y), \mathrm{vec}(X)\rangle$, the characteristic function of a matrix subordinator can be represented as

$$E(\mathrm{e}^{\mathrm{i}\,\mathrm{tr}(L_t^* Z)}) = \exp(t\psi_L(Z)), \qquad Z \in \mathbb{S}_d,$$
$$\text{where } \psi_L(Z) := \mathrm{i}\,\mathrm{tr}(\gamma_L Z) + \int_{\mathbb{S}_d^+} (\mathrm{e}^{\mathrm{i}\,\mathrm{tr}(XZ)} - 1)\nu_L(\mathrm{d}X)$$

with drift $\gamma_L \in \mathbb{S}_d^+$ and Lévy measure $\nu_L$.

The discontinuous part of the quadratic variation of any Lévy process $L$ in $\mathbb{R}^d$,

$$[L, L]_t^{\mathfrak{d}} := \int_0^t \int_{\mathbb{R}^d} xx^* \mu_L(\mathrm{d}s, \mathrm{d}x) = \sum_{0\leq s\leq t} (\Delta L_s)(\Delta L_s)^*$$

is a matrix subordinator with drift zero and Lévy measure given by

$$\nu_{[L,L]^{\mathfrak{d}}}(B) = \int_{\mathbb{R}^d} I_B(xx^*)\nu_L(\mathrm{d}x)$$

for all Borel sets $B \subseteq \mathbb{S}_d$.

## 3. Definition of multivariate COGARCH$(1,1)$ processes

The main idea for the definition of a multivariate COGARCH$(1,1)$ process is to replace the noise $\varepsilon$ of a multivariate GARCH$(1,1)$ process (see [6, 16] and references therein) by the jumps of a multivariate Lévy process $L$, and the autoregressive structure of the covariance matrix process by a continuous-time autoregressive (AR) structure (that is, an Ornstein–Uhlenbeck (OU) type structure). So, the idea is again basically the same as in [10] for the univariate COGARCH$(p,q)$ process.

In the simplest BEKK GARCH$(1,1)$ model of [16], the volatility process is given by

$$\Sigma_n = C + A\Sigma_{n-1}^{1/2}\varepsilon_{n-1}\varepsilon_{n-1}^*\Sigma_{n-1}^{1/2}A^* + B\Sigma_{n-1}B^* \tag{3.1}$$

with $C \in \mathbb{S}_d^+$, $A, B \in M_d(\mathbb{R})$ and $(\varepsilon_n)_{n\in\mathbb{N}_0}$ being an i.i.d. sequence in $\mathbb{R}^d$. This shows that the dynamics of $(\Sigma_n)_{n\in\mathbb{N}_0}$ are those of a multivariate AR process, which is "self-exciting" in the sense that we have an AR structure with the noise given by $(\Sigma_{n-1}^{1/2}\varepsilon_{n-1}\varepsilon_{n-1}^*\Sigma_{n-1}^{1/2})_{n\in\mathbb{N}}$.

Replacing the AR structure with an OU-type structure, using $V_{t-}^{1/2}\,\mathrm{d}[L,L]_t^{\mathfrak{d}}\,V_{t-}^{1/2}$ as "noise", where $L$ is a $d$-dimensional Lévy process, and using the same linear operators as for positive semidefinite processes of OU type (see [5, 35]) now leads to a multivariate continuous-time GARCH$(1,1)$ process $G$ (referred to as a MUCOGARCH$(1,1)$ process in the following) given by the following definition.

***Definition 3.1* (MUCOGARCH(1,1)).** *Let $L$ be an $\mathbb{R}^d$-valued Lévy process and $A, B \in M_d(\mathbb{R})$, $C \in \mathbb{S}_d^{++}$. The process $G = (G_t)_{t\in\mathbb{R}^+}$ solving*

$$\mathrm{d}G_t = V_{t-}^{1/2}\,\mathrm{d}L_t, \tag{3.2}$$

$$V_t = C + Y_t, \tag{3.3}$$

$$\mathrm{d}Y_t = (BY_{t-} + Y_{t-}B^*)\,\mathrm{d}t + AV_{t-}^{1/2}\,\mathrm{d}[L,L]_t^{\mathfrak{d}}\,V_{t-}^{1/2}A^*, \tag{3.4}$$

*with initial values $G_0$ in $\mathbb{R}^d$ and $Y_0$ in $\mathbb{S}_d^+$, is then called a* MUCOGARCH$(1,1)$ process.

*The process $Y = (Y_t)_{t\in\mathbb{R}^+}$ with paths in $\mathbb{S}_d^+$ is referred to as a* MUCOGARCH$(1,1)$ volatility process.

As we are only dealing with MUCOGARCH processes of order $(1,1)$, we often write only "MUCOGARCH" instead of "MUCOGARCH$(1,1)$" in the sequel.

We can also directly state a stochastic differential equation (SDE) for the covariance matrix process $V$:

$$\mathrm{d}V_t = (B(V_{t-} - C) + (V_{t-} - C)B^*)\,\mathrm{d}t + AV_{t-}^{1/2}\,\mathrm{d}[L,L]_t^{\mathfrak{d}}\,V_{t-}^{1/2}A^*. \tag{3.5}$$

This SDE has a "mean-reverting structure" (provided $\sigma(B) \subset (-\infty, 0) + \mathrm{i}\mathbb{R}$), namely, $V$ returns to the level $C$ at an exponential rate determined by $B$, as long as there are no jumps. However, since all jumps are positive semidefinite, as we shall see, $C$ is not a "mean" level, but a lower bound.

Equivalently, we can use the following representation using the vec operator:

$$\mathrm{d}G_t = V_{t-}^{1/2}\,\mathrm{d}L_t, \qquad V_t = C + Y_t,$$

$$\mathrm{d}\,\mathrm{vec}(Y_t) = (B \otimes I + I \otimes B)\,\mathrm{vec}(Y_{t-})\,\mathrm{d}t + (A \otimes A)(V_{t-}^{1/2} \otimes V_{t-}^{1/2})\,\mathrm{d}\,\mathrm{vec}([L,L]_t^{\mathfrak{d}}),$$

$$\mathrm{d}\,\mathrm{vec}(V_t) = (B \otimes I + I \otimes B)(\mathrm{vec}(V_{t-}) - \mathrm{vec}(C))\,\mathrm{d}t + (A \otimes A)(V_{t-}^{1/2} \otimes V_{t-}^{1/2})\,\mathrm{d}\,\mathrm{vec}([L,L]_t^{\mathfrak{d}}).$$

For the MUCOGARCH process to be well defined, it is necessary that there exists a unique solution to the above system of stochastic differential equations and that $V$ does not leave the set $\mathbb{S}_d^+$. In the following, it is implicitly understood that our processes and stochastic differential equations are not living on the space $M_d(\mathbb{R})$ (resp., $\mathbb{R}^{d^2}$), but on the linear subspace $\mathbb{S}_d$ of symmetric matrices (resp., $\mathrm{vec}(\mathbb{S}_d)$). The latter can, as usual, be identified with $\mathbb{R}^{d(d+1)/2}$, when appropriate. The importance of this lies in the fact that $\mathbb{S}_d^{++}$ is an open subset of $\mathbb{S}_d$ and it is most natural to consider (stochastic) differential equations on open sets.

**Theorem 3.2.** *Let $A, B \in M_d(\mathbb{R})$, $C \in \mathbb{S}_d^{++}$ and $L$ be a $d$-dimensional Lévy process. The SDE (3.4) with initial value $Y_0 \in \mathbb{S}_d^+$ then has a unique positive semidefinite solution $(Y_t)_{t\in\mathbb{R}^+}$. The solution $(Y_t)_{t\in\mathbb{R}^+}$ is locally bounded and of finite variation. Moreover, it satisfies $Y_t \geq \mathrm{e}^{Bt} Y_0 \mathrm{e}^{B^*t}$ for all $t \in \mathbb{R}^+$.*

***Remark 3.3.*** (i) An analogous result holds when considering $C \in \mathbb{S}_d^+$ and restricting the initial value to $Y_0 \in \mathbb{S}_d^{++}$ ("locally bounded" needs to be replaced with "locally bounded within $\mathbb{S}_d^+$", as defined in [5], Definition 3.1). All of the following results (except those regarding the existence of stationary solutions) can be immediately adapted to this case.

(ii) For $d = 1$, it is straightforward to see that our definition agrees (after a reparametrization) with the case $p = q = 1$ of the general COGARCH$(p, q)$ definition given in [10], that is, $Y_t$ agrees with their process $\alpha_1 \mathbf{Y}_t$ and our $V_t$ with their $V_{t+}$. Hence, [10], Theorem 2.2, implies that our definition agrees with the original definition given in [26].

Likewise, we could have considered the SDE (3.5). Using the relationship between (3.5) and (3.4), we obtain the following.

**Corollary 3.4.** *Let $A, B \in M_d(\mathbb{R})$, $C \in \mathbb{S}_d^{++}$ and $L$ be a $d$-dimensional Lévy process. Assume that the initial value satisfies $V_0 \geq C$. The SDE (3.5) then has a unique positive definite solution $(V_t)_{t\in\mathbb{R}^+}$ and $V_t \geq C + \mathrm{e}^{Bt}(V_0 - C)\mathrm{e}^{B^*t} \geq C$ for all $t \in \mathbb{R}^+$.*

It may appear natural also to allow initial values $V_0 \in \mathbb{S}_d^+$ with $V_0 < C$. In this case, one still has that $V_t \geq C + \mathrm{e}^{Bt}(V_0 - C)\mathrm{e}^{B^*t}$ needs to be true for any solution of (3.5), as long as it exists. However, in this case, the solution of (3.5) may leave the set $\mathbb{S}_d^+$ and thus have only a finite lifetime, as the following example shows. Take

$$C = \begin{pmatrix} 2 & 0 \\ 0 & 2 \end{pmatrix}, \qquad V_0 = \begin{pmatrix} 0.5 & 0 \\ 0 & 0.5 \end{pmatrix},$$
$$B = \begin{pmatrix} -0.5\ln(10/9) & 0 \\ 1 & -0.5\ln(10/9) \end{pmatrix}, \qquad x = \begin{pmatrix} 1 \\ 1 \end{pmatrix}$$

and $L$ as the zero Lévy process. We then obtain that

$$\mathrm{e}^B = \sqrt{\frac{9}{10}} \begin{pmatrix} 1 & 0 \\ 1 & 1 \end{pmatrix} \quad \text{and} \quad x^* V_1 x = -\frac{11}{4}.$$

So, $V_1 \notin \mathbb{S}_d^+$, although $V_0 \in \mathbb{S}_d^+$. Note that this problem also arises with positive probability if the driving Lévy process is compound Poisson, as it may then happen that there is no jump until time 1.

***Remark 3.5.*** The insight gained from positive semidefinite OU processes in [5] suggests that all eigenvalues of $B$ should have negative real part if one wants a stationary COGARCH volatility process as covariance matrix process. It is clear that, in this case, $V_t \to C$ as $t \to \infty$ if the Lévy process had no jumps. Thus, in general, the process $V$ tends to $C$, as long as the driving Lévy process does not jump. The above counterexample shows that when the process $V$ is smaller than $C$ (in the ordering of the positive semidefinite matrices), this does not occur in a "straight" manner, whereas in the univariate model, the volatility process is always increasing below $C$ (see [27], Proposition 2).

Similar to the usual shot noise representation of OU-type processes, we have the following.

**Theorem 3.6.** *The* MUCOGARCH(1,1) *volatility process* $Y$ *satisfies*

$$Y_t = \mathrm{e}^{Bt} Y_0 \mathrm{e}^{B^* t} + \int_0^t \mathrm{e}^{B(t-s)} A (C + Y_{s-})^{1/2} \, \mathrm{d}[L, L]_s^{\mathfrak{d}} \, (C + Y_{s-})^{1/2} A^* \mathrm{e}^{B^*(t-s)} \tag{3.6}$$

*for all* $t \in \mathbb{R}^+$.

Recently, [38] studied univariate equations of the form $X(t) = J(t) + \int_0^t g(t - s) f(X_{s-}) \, \mathrm{d}Z_s$ and their relation to certain SDEs. In particular, they obtained uniqueness of the solutions under uniform Lipschitz assumptions on $f$. Our equation (3.6) is a multivariate equation of this type, with $f$ being only locally Lipschitz. From the arguments given in [38], one sees that their Theorem 5.2 remains valid in a multivariate setting. Using a localization procedure as in the proof of [40], Theorem 6.6.3, this uniqueness result extends to $f$ being defined only on an open subset and locally Lipschitz. Hence, (3.6) provides an alternative characterization for the MUCOGARCH volatility process.

So far, we have excluded the MUCOGARCH process $G$ itself from the analysis. However, the following result is obtained along the same lines as Theorem 3.2.

**Theorem 3.7.** *Let $A, B \in M_d(\mathbb{R})$, $C \in \mathbb{S}_d^{++}$ and $L$ be a $d$-dimensional Lévy process. The system of SDEs (3.2), (3.4) then has a unique solution $(G_t, Y_t)_{t\in\mathbb{R}^+}$ with paths in $\mathbb{R}^d \times \mathbb{S}_d^+$ for any initial value $(G_0, Y_0)$ in $\mathbb{R}^d \times \mathbb{S}_d^+$. The solution $(G_t, Y_t)_{t\in\mathbb{R}^+}$ is a semimartingale and locally bounded.*

# 4. Properties of the volatility process

## 4.1. Univariate COGARCH(1, 1) bounds

We now show that, similarly to the COGARCH$(p, q)$ case (see [10], Lemma 9.1), the norm of a MUCOGARCH(1, 1) volatility process can be bounded by a univariate COGARCH(1, 1) volatility process. This immediately gives useful conditions for the finiteness of moments and has far-reaching implications regarding the existence of stationary distributions.

In the following, we shall consider a special norm that fits our model particularly well. $\|\cdot\|_2$ denotes the operator norm on $M_{d^2}(\mathbb{R})$ associated with the usual Euclidean norm. Assume, now, that $B$ is diagonalizable and let $S \in GL_d(\mathbb{C})$ be such that $S^{-1}BS$ is diagonal. We then define the norm $\|\cdot\|_{B,S}$ on $M_{d^2}(\mathbb{R})$ by $\|X\|_{B,S} := \|(S^{-1}\otimes S^{-1})X(S\otimes S)\|_2$ for $X \in M_{d^2}(\mathbb{R})$. It should be noted that $\|\cdot\|_{B,S}$ depends both on $B$ and on the choice of the matrix $S$ diagonalizing $B$. Actually, $\|\cdot\|_{B,S}$ is again an operator norm, namely the one associated with the norm $\|x\|_{B,S} := \|(S^{-1}\otimes S^{-1})x\|_2$ on $\mathbb{R}^{d^2}$. Besides, $\|\cdot\|_{B,S}$ is simply the norm $\|\cdot\|_2$, provided $S$ is a unitary matrix (see [23], page 308).

**Theorem 4.1.** *Let $Y$ be a* MUCOGARCH *volatility process with initial value $Y_0 \in \mathbb{S}_d^+$ and driven by a Lévy process $L$ in $\mathbb{R}^d$. Assume, further, that $B \in M_d(\mathbb{R})$ is diagonalizable and let $S \in GL_d(\mathbb{C})$ be such that $S^{-1}BS$ is diagonal. The process solving the SDE,*

$$\begin{aligned} \mathrm{d}y_t &= 2\lambda y_{t-}\,\mathrm{d}t + \|S\|_2^2\|S^{-1}\|_2^2 K_{2,B}\|A\otimes A\|_{B,S}\left(\frac{\|C\|_2}{K_{2,B}} + y_{t-}\right)\mathrm{d}\tilde{L}_t, \\ y_0 &= \|\operatorname{vec}(Y_0)\|_{B,S} \end{aligned} \tag{4.1}$$

*with*

$$\tilde{L}_t := \int_0^t \int_{\mathbb{R}^d} \|\operatorname{vec}(xx^*)\|_{B,S}\,\mu_L(\mathrm{d}s, \mathrm{d}x), \qquad \lambda := \max(\Re(\sigma(B)))$$

*and*

$$K_{2,B} := \max_{X\in\mathbb{S}_d^+, \|X\|_2=1}\left(\frac{\|X\|_2}{\|\operatorname{vec}(X)\|_{B,S}}\right),$$

*is the volatility process of a univariate* MUCOGARCH(1, 1) *process and $y$ satisfies*

$$\|\operatorname{vec}(Y_t)\|_{B,S} \le y_t \qquad \text{for all } t\in\mathbb{R}^+ \text{ a.s.} \tag{4.2}$$

*Moreover, $K_{2,B} \le \|S\|_2^2 \max_{X\in\mathbb{S}_d^+, \|X\|_2=1}(\frac{\|X\|_2}{\|\operatorname{vec}(X)\|_2}) \le \|S\|_2^2$.*

***Remark 4.2.*** (i) Provided $S$ is unitary, $K_{2,B} = 1$. Otherwise, an inspection of the proof shows that the inequality (4.2) also holds if $K_{2,B}$ is replaced by $\|S\|_2$ in (4.1), which saves one from calculating the value of $K_{2,B}$ in practice. Likewise, $\|A \otimes A\|_{B,S}$ can be replaced by $\|A \otimes A\|_2 = \|A\|_2^2$ since $\|(A \otimes A)((C + Y_{\Gamma_1 -})^{1/2} \otimes (C + Y_{\Gamma_1 -})^{1/2})\|_{B,S} \leq \|S\|_2^2 \|S^{-1}\|_2^2 \|(A \otimes A)((C + Y_{\Gamma_1 -})^{1/2} \otimes (C + Y_{\Gamma_1 -})^{1/2})\|_2$.

This can be done in all of the forthcoming results involving $K_{2,B}$ or $\|A \otimes A\|_{B,S}$.

(ii) As can be seen from the proof, the diagonalizability of $B$ is essential and, unfortunately, it seems very intricate to extend the result to the non-diagonalizable case. In applications, however, this appears to be no severe constraint, as the non-diagonalizable matrices have Lebesgue measure zero.

Since the finiteness of moments of univariate COGARCH$(1, 1)$ processes is well known from [26], Section 4, we can now give sufficient conditions for the MUCOGARCH volatility process to have some finite moments, which we will improve upon in Proposition 4.7.

**Proposition 4.3.** *Let $k \in \mathbb{N}$, $Y_0 \in \mathbb{S}_d^+$ such that $E(\|Y_0\|^k) < \infty$ and let $B$ be diagonalizable. Assume further that the* MUCOGARCH *volatility process $Y$ is driven by a Lévy process $L$ satisfying $E(\|L_1\|^{2k}) < \infty$. Then $E(\|Y_t\|^k) < \infty$ for all $t \in \mathbb{R}^+$ and $t \mapsto E(\|Y_t\|^k)$ is locally bounded.*

## 4.2. Markovian properties and stationarity

Turning to the study of the Markovian properties of a MUCOGARCH process, we refer to standard references like [15, 17, 18] for definitions and necessary general results. Moreover, we implicitly assume that our given filtered probability space is enlarged as in [36], page 293, to allow for arbitrary initial conditions of the SDEs, and the weak Feller property is defined as in [15], namely by demanding that the transition semigroup is stochastically continuous and maps the bounded continuous functions on the state space into themselves.

The usual results on the Markov properties of SDEs (see [36], Section V.6) extend to locally Lipschitz SDEs on open sets (see [40], Section 6.7.1.2, for details) and to closed sets, provided the solution is ensured to stay in the closed set at all times and the SDE is defined on an open set containing the closed set. The latter is the case for the MUCOGARCH, the closed set $\mathbb{S}_d^+$ and the open set $U_{C,\varepsilon}$ as defined in the proof of Theorem 3.2. Thus, we obtain the following result.

**Theorem 4.4.** *The* MUCOGARCH *process $(G, Y)$ and its volatility process $Y$ alone are temporally homogeneous strong Markov processes on $\mathbb{R}^d \times \mathbb{S}_d^+$ and $\mathbb{S}_d^+$, respectively, and they have the weak Feller property.*

One of the most important questions regarding Markov processes in applications is the existence of stationary distributions.

**Theorem 4.5.** *Let $B \in M_d(\mathbb{R})$ be diagonalizable with $S \in GL_d(\mathbb{C})$ such that $S^{-1}BS$ is diagonal. Furthermore, let $L$ be a $d$-dimensional Lévy process with non-zero Lévy measure, $\lambda$ be defined as in Theorem 4.1 and $\alpha_1 := \|S\|_2^2 \|S^{-1}\|_2^2 K_{2,B} \|A \otimes A\|_{B,S}$. Assume that*

$$\int_{\mathbb{R}^d} \log(1 + \alpha_1 \| \operatorname{vec}(yy^*)\|_{B,S}) \nu_L(\mathrm{d}y) < -2\lambda. \tag{4.3}$$

*There then exists a stationary distribution $\mu \in \mathcal{M}_1(\mathbb{S}_d^+)$, that is, the set of all probability measures on the Borel-$\sigma$-algebra of $\mathbb{S}_d^+$, for the* MUCOGARCH(1, 1) *volatility process $Y$ such that*

$$\int_{\mathbb{R}^d} ((1 + \alpha_1 \| \operatorname{vec}(yy^*)\|_{B,S})^k - 1) \nu_L(\mathrm{d}y) < -2\lambda k \tag{4.4}$$

*for some $k \in \mathbb{N}$ implies that $\int_{\mathbb{S}_d^+} \|x\|^k \mu(\mathrm{d}x) < \infty$, that is, that the $k$th moment of the stationary distribution is finite.*

Of course, this result immediately translates to stationarity of $V$. Moreover, for $d = 1$, it recovers the necessary and sufficient stationarity condition of [26].

***Remark 4.6.*** (a) From [26], Lemma 4.1(d), it follows that, if (4.4) is satisfied for $k \in \mathbb{N}$, then it is also satisfied for all $\tilde{k} \in \mathbb{N}$, $\tilde{k} \leq k$.

(b) Combining the results of Section 4.1 shows that $\alpha_1 = \|A\|_2^2$ and $\| \cdot \|_{B,S} = \| \cdot \|_2$ if $B$ is normal.

Establishing uniqueness of the stationary distribution and convergence to the stationary distribution for arbitrary starting values appears to be a rather intricate question due to the Lipschitz property holding only locally and the fact that $\mathrm{d}[L, L]^{\mathfrak{d}}$ lives on the rank one matrices. However, in the next section, we obtain at least asymptotic second-order stationarity and that the stationary second-order structure is unique under some technical conditions.

To conclude this section, we consider some examples exploring the dependence of the stationary distribution on the parameters and the relation to the stationarity of univariate COGARCH processes.

***Example 4.1.*** Let $c \in \mathbb{R}^+ \setminus \{0\}$, $A, B \in M_d(\mathbb{R})$ and $L$ be a $d$-dimensional Lévy process. If $V$ satisfies

$$\mathrm{d}V_t = (B(V_{t-} - cI_d) + (V_{t-} - cI_d)B^*)\,\mathrm{d}t + AV_{t-}^{1/2}\,\mathrm{d}[L, L]_t^{\mathfrak{d}}\, V_{t-}^{1/2} A^*, \tag{4.5}$$

then $Z$ defined by $Z_t = V_t/c$ satisfies

$$\mathrm{d}Z_t = (B(Z_{t-} - I_d) + (Z_{t-} - I_d)B^*)\,\mathrm{d}t + AZ_{t-}^{1/2}\,\mathrm{d}[L, L]_t^{\mathfrak{d}}\, Z_{t-}^{1/2} A^*, \tag{4.6}$$

which does not depend on $c$. In particular, if $\mu \in \mathcal{M}_1(\mathbb{S}_d^+)$ is a stationary distribution for (4.6), then $\mu_c \in \mathcal{M}_1(\mathbb{S}_d^+)$, defined by $\mu_c(W) = \mu(W/c)$ for all Borel sets $W \subset \mathbb{S}_d^+$, is a stationary distribution for (4.5).

***Example 4.2.*** Assume that $A, B, C$ and $Y_0$ are diagonal, the components of $Y_0$ are independent and the components of the Lévy process are completely independent, that is, whenever $L$ has a jump, then only one of the $d$ components jumps. In this case, $Y$ or $V$, respectively, consists of $d$ independent univariate COGARCH$(1, 1)$ volatility processes. If each of the $d$ univariate COGARCH$(1, 1)$ volatility processes converges in distribution to a stationary distribution, then $Y$ or $V$, respectively, converges in distribution to a stationary distribution. In this example, condition (4.3) can be shown to imply the necessary and sufficient stationarity condition of [26], Theorem 3.1, for all components simultaneously. Actually, condition (4.3) is stronger than requiring that the univariate stationarity condition be satisfied for all components.

However, it should be noted that the picture is very different when $Y_0$ is not diagonal because then jumps in one component of $L$ typically affect all components of $Y$. Hence, it is not clear whether one still has convergence to a stationary distribution and whether this has to be the same distribution as the limit distribution when $Y_0$ is diagonal. When we have that $Y$ is asymptotically second order stationary (see Theorem 4.20 below) and the limiting distribution for a diagonal $Y_0$ has finite second moments, the off-diagonal (covariance) elements of $Y$ or $V$, respectively, necessarily converge to zero in $L^2$ as $t \to \infty$.

***Example 4.3.*** A degenerate situation occurs if $C \in \mathbb{S}_d^+ \backslash \mathbb{S}_d^{++}$. Take $d = 2$, $A = \alpha I_2$, $B = -\beta I_2$ with $\alpha, \beta \in \mathbb{R}^+ \backslash \{0\}$ and $C = \left(\begin{smallmatrix} 1 & 1 \\ 1 & 1 \end{smallmatrix}\right)$. If one has that $Y_0 = yC$ with $y \in \mathbb{R}^+$ (possibly random), then $Y$ is at all times a scalar multiple of $C$ and, when $L$ jumps, all components (variance and covariance ones) have a jump of the same height. However, one again has a completely different picture if $Y_0$ is chosen differently, for example, $Y$ is in $\mathbb{S}_d^{++}$ at all times, provided $Y_0 \in \mathbb{S}_d^{++}$. This also shows that the assumption $C \in \mathbb{S}_d^{++}$ made in the definition of the MUCOGARCH processes is essential to avoid pathological situations.

## 4.3. Second-order moment structure

Assuming stationarity and the existence of the relevant moments of the stationary solution, we calculate explicit expressions for the moments of a stationary MUCOGARCH$(1, 1)$ volatility process in this section, treat the non-stationary case along the way and present results regarding (asymptotic) second-order stationarity. Due to the special structure of the stochastic differential equation (3.4), especially due to the presence of the matrix square root, it is only possible under certain assumptions on the Lévy process to obtain explicit formulae. The results of this chapter provide the basis for method of moments estimation, provided the volatility process is (approximately) observed and shows that the second-order structure of the volatility process is in line with observed financial data, since the matrix exponential decay of the autocovariance is rather flexible and has been found realistic in the analysis of OU-type models (see [34]).

Henceforth, we often assume the following in this section.

***Assumption 4.1.*** $(Y_t)_{t\in\mathbb{R}^+}$ *is a second-order stationary* MUCOGARCH$(1, 1)$ *volatility process.*

***Assumption 4.2.*** *The pure jump part of the driving Lévy process $(L_t)_{t\in\mathbb{R}^+}$ has finite variance which is a scalar multiple of the identity matrix:*

*If we let $L_t^{\mathfrak{d}} := \int_0^t \int_{\|x\|\le 1} x(\mu_L(\mathrm{d}s,\mathrm{d}x) - \mathrm{d}s\,\nu_L(\mathrm{d}x)) + \int_0^t \int_{\|x\|>1} x\mu_L(\mathrm{d}s,\mathrm{d}x)$ denote the pure jump part of $L$, then this means that there exists a $\sigma_L \in \mathbb{R}^+$ such that $\mathrm{var}(L_1^{\mathfrak{d}}) = \int_{\mathbb{R}^d} xx^* \nu_L(\mathrm{d}x) = \sigma_L I_d$.*

This assumption is comparable to considering only *standard* Brownian motion in Brownian-motion-based models and, hence, not very restrictive since any Lévy process with finite second moments can be transformed into one satisfying Assumption 4.2 by a linear transformation and since the variance of $G$ can still be flexibly modelled via the remaining parameters, as will be seen from Proposition 5.2.

First, we need a refinement of Proposition 4.3 to the case where $B$ is not diagonalizable.

**Proposition 4.7.** *Let $Y$ be a* MUCOGARCH(1,1) *volatility process and $k \in \{1\} \cup [2,\infty)$.*

*If $E(\|Y_0\|^k) < \infty$ and $E(\|L_1\|^{2k}) < \infty$, then $E(\|Y_t\|^k) < \infty$ for all $t \in \mathbb{R}^+$ and $t \mapsto E(\|Y_t\|^k)$ is locally bounded.*

We can now calculate the expected value of the volatility.

**Theorem 4.8.** *Assume that Assumption 4.2 holds:*

(i) *If $E(\|Y_0\|) < \infty$, then*

$$E(\mathrm{vec}(Y_t)) = \mathrm{e}^{\mathscr{B}t} E(\mathrm{vec}(Y_0)) + \int_0^t \mathrm{e}^{\mathscr{B}(t-s)}\,\mathrm{d}s\,\sigma_L(A\otimes A)\,\mathrm{vec}(C)$$

*with $\mathscr{B} := B\otimes I_d + I_d \otimes B + \sigma_L A\otimes A$. If $\mathscr{B}$ is invertible, then*

$$E(\mathrm{vec}(Y_t)) = \mathrm{e}^{\mathscr{B}t}(E(\mathrm{vec}(Y_0)) + \sigma_L \mathscr{B}^{-1}(A\otimes A)\,\mathrm{vec}(C)) - \sigma_L\mathscr{B}^{-1}(A\otimes A)\,\mathrm{vec}(C) \quad (4.7)$$

*for all $t \in \mathbb{R}^+$.*

(ii) *Under Assumption 4.1, the stationary expected value $E(Y_0)$ of the* MUCOGARCH *volatility process satisfies*

$$BE(Y_0) + E(Y_0)B^* + \sigma_L AE(Y_0)A^* = -\sigma_L ACA^*. \quad (4.8)$$

*If $\mathscr{B}$ is invertible, then the following formulae hold:*

$$\begin{aligned} E(\mathrm{vec}(Y_0)) &= -\sigma_L\mathscr{B}^{-1}(A\otimes A)\,\mathrm{vec}(C) \quad \text{and} \\ E(\mathrm{vec}(V_0)) &= \mathscr{B}^{-1}(B\otimes I_d + I_d\otimes B)\,\mathrm{vec}(C). \end{aligned} \quad (4.9)$$

***Remark 4.9.*** Observe that the stationary expectation is the limit of the expected value in (i) for $t\to\infty$ provided $\sigma(\mathscr{B}) \subset (-\infty, 0) + i\mathbb{R}$.

Theorem 4.5 can not only be used to show that Assumption 4.1 is satisfied, but also to ensure the invertibility of $\mathscr{B}$. Hence, Theorem 4.8 provides an explicit expression for the mean of the stationary distribution of Theorem 4.5.

**Lemma 4.10.** *Assume that (4.4) is satisfied with $k=1$ for the* MUCOGARCH *volatility process $Y$ and that Assumption 4.2 holds. Then, $\mathscr{B}$, defined as above, is invertible and $\sigma(\mathscr{B}) \subset (-\infty,0)+\mathrm{i}\mathbb{R}$.*

Before analyzing the variance, let us study the autocovariance function. If $(X_t)_{t\in\mathbb{R}^+}$ is a second-order stationary process with values in $\mathbb{R}^d$, the autocovariance function $\mathrm{acov}_X : \mathbb{R} \mapsto M_d(\mathbb{R})$ of $X$ is given by $\mathrm{acov}_X(h) = \mathrm{cov}(X_h, X_0) = E(X_h X_0^*) - E(X_0)E(X_0)^*$ for $h \geq 0$ and by $\mathrm{acov}_X(h) = (\mathrm{acov}_X(-h))^*$ for $h < 0$. As we are considering matrix-valued processes $(Z_t)_{t\in\mathbb{R}}$ in the following, we set $\mathrm{acov}_Z := \mathrm{acov}_{\mathrm{vec}(Z)}$ in this case.

**Theorem 4.11.** (i) *Under Assumptions 4.1 and 4.2, the autocovariance function of the* MUCOGARCH *volatility process satisfies*

$$\frac{\mathrm{d}}{\mathrm{d}h}\mathrm{acov}_Y(h) = (B\otimes I_d + I_d \otimes B + \sigma_L A\otimes A)\,\mathrm{acov}_Y(h) \tag{4.10}$$

*for $h \geq 0$.*

*Hence,*

$$\mathrm{acov}_Y(h) = \mathrm{acov}_V(h) = \mathrm{e}^{(B\otimes I_d + I_d\otimes B+\sigma_L A\otimes A)h}\,\mathrm{var}(\mathrm{vec}(Y_0)), \qquad h \geq 0. \tag{4.11}$$

(ii) *If Assumption 4.2 is satisfied and $E(\|Y_0\|^2), E(\|L_1\|^4)$ are finite, it holds that*

$$\mathrm{cov}(Y_{u+h}, Y_u) = \mathrm{cov}(V_{u+h}, V_u) = \mathrm{e}^{(B\otimes I_d + I_d\otimes B+\sigma_L A\otimes A)h}\,\mathrm{var}(\mathrm{vec}(Y_u)) \tag{4.12}$$

*for all $u, h \geq 0$.*

The autocovariance function of the volatility process $Y$ is thus exponentially decreasing in a matrix sense, so the individual entries may decay as sums of exponentials, exponentially damped sinusoids (if the eigenvalues have non-vanishing complex parts) or exponentially damped polynomials (if the matrix is not diagonalizable).

However, we are so far lacking an explicit expression for $\mathrm{var}(\mathrm{vec}(Y_0))$. Unfortunately, our Assumption 4.2 on the second moment of the jumps of the driving Lévy process $L$ seems not to be sufficient to obtain an explicit expression for the variance.

As we shall see from the proofs, the quadratic variation of the vectorized discontinuous part of the quadratic variation of the driving Lévy process,

$$[\mathrm{vec}([L,L]^{\mathfrak{d}}), \mathrm{vec}([L,L]^{\mathfrak{d}})]_t^{\mathfrak{d}} = \int_0^t \int_{\mathbb{R}^d} \mathrm{vec}(xx^*)\,\mathrm{vec}(xx^*)^* \mu_L(\mathrm{d}s, \mathrm{d}x),$$

which is again a pure jump Lévy process of finite variation, will appear in our calculations of the second moment and we need it to have finite expectation. In fact, we even need to make specific assumptions on its expectation.

To determine what assumptions are reasonable, let us assume for a moment that $L$ is a $d$-dimensional compound Poisson process with rate one and with jump distribution being the $d$-dimensional standard normal distribution. This implies that $[L, L]^{\mathfrak{d}}$ is a compound Poisson process with rate one and with the jump distribution being a Wishart distribution. Then, denoting the $d$-dimensional standard normal distribution by $N(\mathrm{d}x)$ and noting that $\mathrm{vec}(xx^*)\,\mathrm{vec}(xx^*)^* = (x \otimes x)(x^* \otimes x^*) = (xx^*) \otimes (xx^*)$, we have

$$\begin{aligned} E([\mathrm{vec}([L,L]^{\mathfrak{d}}), \mathrm{vec}([L,L]^{\mathfrak{d}})]_1^{\mathfrak{d}}) &= \int_{\mathbb{R}^d} (xx^*) \otimes (xx^*) N(\mathrm{d}x) \\ &= I_{d^2} + K_d + \mathrm{vec}(I_d)\,\mathrm{vec}(I_d)^*, \end{aligned} \tag{4.13}$$

from [30], Theorem 4.1. Here, $K_d \in M_{d^2}(\mathbb{R})$ denotes the commutation matrix which can be characterized by $K_d\,\mathrm{vec}(A) = \mathrm{vec}(A^*)$ for all $A \in M_d(\mathbb{R})$ (see [30] for more details). This generalizes to the following result.

**Lemma 4.12.** *Let $L$ be a $d$-dimensional compound Poisson process with rate $c$ and with jumps distributed like $\sqrt{\varepsilon}X$, where $X$ is a $d$-dimensional standard normal random variable and $\varepsilon$ is a random variable in $\mathbb{R}^+$ with finite variance and independent of $X$. Then,*

$$E([\mathrm{vec}([L,L]^{\mathfrak{d}}), \mathrm{vec}([L,L])^{\mathfrak{d}}]_1^{\mathfrak{d}}) = cE(\varepsilon^2)(I_{d^2} + K_d + \mathrm{vec}(I_d)\,\mathrm{vec}(I_d)^*). \tag{4.14}$$

Moving away from a Lévy process of finite activity, a similar result holds for the following variant of type $G$ processes, a special kind of a normal mixture.

***Definition 4.13 (Type $\widetilde{G}$).*** *Let $L$ be a $d$-dimensional Lévy process. If there exists an $\mathbb{R}^+$-valued infinitely divisible random variable $\varepsilon$ independent of a $d$-dimensional standard normal random variable $X$ such that $L_1 \stackrel{\mathscr{L}}{=} \sqrt{\varepsilon}X$, then $L$ is said to be of* type $\widetilde{G}$. *(Here, $\stackrel{\mathscr{L}}{=}$ denotes equality in law.)*

We have chosen the name "type $\widetilde{G}$" above because these processes correspond to a particular case of multG laws as defined in [2], Definition 3.1. Actually, many interesting Lévy processes are of type $\widetilde{G}$, for instance, the multivariate symmetric GH (NIG) processes with the parameter $\Sigma$ set to $I_d$ (see [9, 31]). For details on distributions and Lévy processes of type $G$ in general, we refer to [2, 29].

**Lemma 4.14.** *Let $L$ be a $d$-dimensional Lévy process of type $\widetilde{G}$ with a finite fourth moment. Then $E([\mathrm{vec}([L,L]^{\mathfrak{d}}), \mathrm{vec}([L,L]^{\mathfrak{d}})]_1^{\mathfrak{d}}) = \rho_L(I_{d^2} + K_d + \mathrm{vec}(I_d)\,\mathrm{vec}(I_d)^*)$ with $\rho_L \in \mathbb{R}^+$.*

These results motivate the following assumption.

***Assumption 4.3.*** *The pure jump part of the driving Lévy process $(L_t)_{t\in\mathbb{R}^+}$ has a finite fourth moment, that is, $\int_{\mathbb{R}^d}\|x\|^4\nu_L(\mathrm{d}x)<\infty$, and there exists a real constant $\rho_L$ such that*

$$E([\operatorname{vec}([L,L]^{\mathfrak{d}}),\operatorname{vec}([L,L]^{\mathfrak{d}})]_1^{\mathfrak{d}})=\rho_L(I_{d^2}+K_d+\operatorname{vec}(I_d)\operatorname{vec}(I_d)^*).$$

Intuitively, this means that the jumps of $L$ have the same fourth moment as a standard normal distribution. This assumption is considerably more restrictive than Assumption 4.2. However, it is comparable to the assumptions made for discrete-time multivariate GARCH processes (see [21]) and from the proofs, one sees that explicit results are only obtainable if the fourth moment of the jumps is comprised of well-understood matrices which act on tensor products in a suitable way.

To state our next result, we need to introduce some additional special linear operators and matrices. If we define

$$\begin{aligned}&\mathbf{Q}: M_{d^2}(\mathbb{R})\to M_{d^2}(\mathbb{R}),\\ &(\mathbf{Q}X)_{(k-1)d+l,(p-1)d+q}=X_{(k-1)d+p,(l-1)d+q}\qquad \text{for all } k,l,p,q=\{1,2,\ldots,d\},\end{aligned}\tag{4.15}$$

then $\mathbf{Q}^{-1}=\mathbf{Q}$, obviously, and $\mathbf{Q}(\operatorname{vec}(X)\operatorname{vec}(Z)^*)=X\otimes Z$ for all $X,Z\in\mathbb{S}_d$ (see [34], Theorem 4.3). Furthermore, we define $\mathcal{Q}\in M_{d^4}(\mathbb{R})$ as the matrix associated with the linear map $\operatorname{vec}\circ\mathbf{Q}\circ\operatorname{vec}^{-1}$ on $\mathbb{R}^{d^4}$ and $\mathcal{K}_d\in M_{d^4}(\mathbb{R})$ as the matrix associated with the linear map $\operatorname{vec}(K_d\operatorname{vec}^{-1}(x))$ for $x\in\mathbb{R}^{d^4}$, where $\operatorname{vec}: M_{d^2}(\mathbb{R})\to\mathbb{R}^{d^4}$. It is easy to see that both $\mathcal{Q}$ and $\mathcal{K}_d$ simply permute the entries of a vector $x\in\mathbb{R}^{d^4}$. Both $\mathcal{Q}$ and $\mathcal{K}_d$ are thus permutation matrices, so we have $\|\mathcal{Q}\|_2=\|\mathcal{K}_d\|_2=1$, where $\|\cdot\|_2$ is the operator norm associated with the usual Euclidean norm on $\mathbb{R}^{d^4}$.

**Theorem 4.15.** *Assume that Assumptions 4.2 and 4.3 hold.*

(i) *If $E(\|Y_0\|^2)<\infty$, then*

$$\begin{aligned}&\frac{\mathrm{d}}{\mathrm{d}t}\operatorname{vec}(E(\operatorname{vec}(Y_t)\operatorname{vec}(Y_t)^*))\\ &\quad=\frac{\mathrm{d}}{\mathrm{d}t}E(\operatorname{vec}(Y_t)\otimes\operatorname{vec}(Y_t))\\ &\quad=\mathscr{C}\operatorname{vec}(E(\operatorname{vec}(Y_t)\operatorname{vec}(Y_t)^*))+(\sigma_L(A\otimes A)\otimes I_{d^2}+\mathscr{A}\mathcal{R})\operatorname{vec}(C)\otimes E(\operatorname{vec}(Y_t))\\ &\qquad+(\sigma_L I_{d^2}\otimes(A\otimes A)+\mathscr{A}\mathcal{R})E(\operatorname{vec}(Y_t))\otimes\operatorname{vec}(C)+\mathscr{A}\mathcal{R}\operatorname{vec}(C)\otimes\operatorname{vec}(C),\end{aligned}\tag{4.16}$$

*where*

$$\begin{aligned}\mathscr{A}&=(A\otimes A)\otimes(A\otimes A);\qquad \mathcal{R}=\rho_L(\mathcal{Q}+\mathcal{K}_d\mathcal{Q}+I_{d^4})\\ \mathscr{C}&:=(B\otimes I_d+I_d\otimes B)\otimes I_{d^2}+I_{d^2}\otimes(B\otimes I_d+I_d\otimes B)\\ &\quad+\sigma_L((A\otimes A)\otimes I_{d^2}+I_{d^2}\otimes(A\otimes A))+\mathscr{A}\mathcal{R}.\end{aligned}$$

(ii) *Under Assumption 4.1, the stationary second moment* $E(\mathrm{vec}(Y_0)\,\mathrm{vec}(Y_0)^*)$ *of the* MUCOGARCH *volatility process satisfies*

$$\begin{aligned}&\mathscr{B}E(\mathrm{vec}(Y_0)\,\mathrm{vec}(Y_0)^*)+E(\mathrm{vec}(Y_0)\,\mathrm{vec}(Y_0)^*)\mathscr{B}^*\\&\qquad+(A\otimes A)\mathbf{R}E(\mathrm{vec}(Y_0)\,\mathrm{vec}(Y_0)^*)(A^*\otimes A^*)\\&\quad=-\sigma_L[(A\otimes A)\,\mathrm{vec}(C)E(\mathrm{vec}(Y_0))^*+E(\mathrm{vec}(Y_0))\,\mathrm{vec}(C)^*(A^*\otimes A^*)]\\&\qquad-(A\otimes A)\mathbf{R}(E(\mathrm{vec}(Y_0))\,\mathrm{vec}(C)^*+\mathrm{vec}(C)E(\mathrm{vec}(Y_0))^*+\mathrm{vec}(C)\,\mathrm{vec}(C)^*)(A^*\otimes A^*),\end{aligned}\tag{4.17}$$

*with* $\mathbf{R}:=\rho_L(\mathbf{Q}+K_d\mathbf{Q}+I_{d^2})$.

*Provided* $\mathscr{C}$ *is invertible,* $E(\mathrm{vec}(Y_0)\,\mathrm{vec}(Y_0)^*)$ *is given by*

$$\begin{aligned}&\mathrm{vec}(E(\mathrm{vec}(Y_0)\,\mathrm{vec}(Y_0)^*))\\&\quad=-\mathscr{C}^{-1}[\mathscr{A}\mathcal{R}(\mathrm{vec}(C)\otimes\mathrm{vec}(C))\\&\qquad\qquad+(\sigma_L(A\otimes A)\otimes I_{d^2}+\mathscr{A}\mathcal{R})\,\mathrm{vec}(C)\otimes E(\mathrm{vec}(Y_0))\\&\qquad\qquad+(\sigma_L I_{d^2}\otimes(A\otimes A)+\mathscr{A}\mathcal{R})E(\mathrm{vec}(Y_0))\otimes\mathrm{vec}(C)].\end{aligned}\tag{4.18}$$

***Remark 4.16.*** The differential equation (4.16) is an inhomogeneous linear differential equation with constant coefficients. Hence, it is standard to obtain an explicit solution. We refrain from stating it, as the stationary case seems to be of the most importance.

Again, condition (4.4) of Theorem 4.5, which ensures the existence of moments of the stationary distribution obtained there, also implies invertibility of $\mathscr{C}$ under an additional technical assumption.

To state the result, we set $\mathcal{S}=S\otimes S\otimes S\otimes S$ and define a new norm $\|\cdot\|_{\widetilde{B,S}}$ on $\mathbb{R}^{d^4}$ by setting $\|x\|_{\widetilde{B,S}}=\|\mathcal{S}^{-1}x\|_2$. The associated operator norm on $M_{d^4}(\mathbb{R})$ is given by $\|X\|_{\widetilde{B,S}}=\|\mathcal{S}^{-1}X\mathcal{S}\|_2$.

**Lemma 4.17.** *Assume that (4.4) is satisfied with* $k=2$ *for the* MUCOGARCH *volatility process* $Y$ *and that Assumptions 4.2 and 4.3 hold. Provided that*

$$\|\mathcal{Q}+\mathcal{K}_d\mathcal{Q}+I_{d^4}\|_{\widetilde{B,S}}\le K_{2,B}^2\|\mathrm{vec}(I_{d^2}+K_d+\mathrm{vec}(I_d)\,\mathrm{vec}(I_d)^*)\|_{\widetilde{B,S}}\tag{4.19}$$

*also holds,* $\sigma(\mathscr{C})\subset(-\infty,0)+\mathrm{i}\mathbb{R}$ *and* $\mathscr{C}$ *is invertible.*

A rather unpleasant feature of this lemma is that we need the technical condition (4.19). The following lemma shows that it is always true if $S$ is unitary, and for concrete parameter values, it can, of course, be checked numerically.

**Lemma 4.18.** *If* $S$ *is unitary, then (4.19) is satisfied.*

We end our comprehensive calculations for the second-order moment structure of the MUCOGARCH volatility process by turning to the stationary variance.

**Corollary 4.19.** *If Assumptions 4.1–4.3 hold and $\mathscr{C}$ is invertible, then the stationary variance* $\mathrm{var}(\mathrm{vec}(Y_0)) = \mathrm{var}(\mathrm{vec}(V_0))$ *of the* MUCOGARCH *volatility process is given by*

$$\begin{aligned} \mathrm{vec}(\mathrm{var}(\mathrm{vec}(Y_0))) = -\mathscr{C}^{-1}[&(\sigma_L^2 \mathscr{C}(\mathscr{B}^{-1} \otimes \mathscr{B}^{-1})\mathscr{A} + \mathscr{A}\mathcal{R})(\mathrm{vec}(C) \otimes \mathrm{vec}(C)) \\ &+ (\sigma_L (A \otimes A) \otimes I_{d^2} + \mathscr{A}\mathcal{R})\,\mathrm{vec}(C) \otimes E(\mathrm{vec}(Y_0)) \\ &+ (\sigma_L I_{d^2} \otimes (A \otimes A) + \mathscr{A}\mathcal{R}) E(\mathrm{vec}(Y_0)) \otimes \mathrm{vec}(C)]. \end{aligned} \tag{4.20}$$

**Proof.** Combine (4.18),

$$\mathrm{vec}(E(\mathrm{vec}(Y_0))E(\mathrm{vec}(Y_0)^*)) = \sigma_L^2 (\mathscr{B}^{-1} \otimes \mathscr{B}^{-1})\mathscr{A}(\mathrm{vec}(C) \otimes \mathrm{vec}(C))$$

and the elementary formula $\mathrm{var}(\mathrm{vec}(Y_0)) = E(\mathrm{vec}(Y_0)\,\mathrm{vec}(Y_0)^*) - E(\mathrm{vec}(Y_0))E(\mathrm{vec}(Y_0)^*)$. □

Under specific moment assumptions on the driving Lévy process, we have thus calculated the second-order structure of a stationary MUCOGARCH volatility process completely.

Finally, we give conditions ensuring (asymptotic) second-order stationarity. A stochastic process $X$ in $\mathbb{S}_d$ is said to be asymptotically second-order stationary with mean $\mu \in \mathbb{R}^{d^2}$, variance $\Sigma \in \mathbb{S}_{d^2}^+$ and autocovariance function $f : \mathbb{R}^+ \to M_{d^2}(\mathbb{R})$ if it has finite second moments and

$$\lim_{t\to\infty} E(X_t) = \mu, \qquad \lim_{t\to\infty} \mathrm{var}(\mathrm{vec}(X_t)) = \Sigma \quad \text{and}$$
$$\lim_{t\to\infty} \sup_{h\in\mathbb{R}^+} \{\| \mathrm{cov}(\mathrm{vec}(X_{t+h}), \mathrm{vec}(X_t)) - f(h)\|\} = 0.$$

**Theorem 4.20.** *Let Assumptions 4.2 and 4.3 be satisfied and assume further that the matrices $B, \mathscr{B}, \mathscr{C}$ are such that $\sigma(B), \sigma(\mathscr{B}), \sigma(\mathscr{C}) \subset (-\infty, 0) + i\mathbb{R}$.*

(i) *If $Y_0$ satisfies (4.9) and (4.20), then the* MUCOGARCH *volatility process $Y$ is second-order stationary.*

(ii) *If $E(\|Y_0\|^2) < \infty$, then the* MUCOGARCH *volatility process $Y$ is asymptotically second-order stationary with mean, variance and autocovariance function given by (4.9), (4.20) and (4.11).*

# 5. The increments of the MUCOGARCH(1, 1) process

Thus far, we have mainly studied the MUCOGARCH volatility processes $Y$ and $V$. However, in practice, one typically cannot observe the volatility, but only the process $G$ (which, in a financial context, for instance, resembles log-prices) at finitely many points in time. In the following, we presume that $G$ is observed on a discrete-time grid starting at zero and with fixed grid size $\Delta > 0$. It is obvious how the results of this section generalize to non-equidistant observations or to the setup considered in [10, 26].

In financial time series, one commonly observes that the returns themselves are uncorrelated, but the "squared returns" (that is, the return vector times its transpose in a multivariate setting) are considerably correlated. The following results show that the MUCOGARCH model can reproduce this very important stylized feature and, furthermore, they provide the basis for simple moment estimators (as in [22]).

We define the sequence of increments $\mathbf{G} = (\mathbf{G}_n)_{n\in\mathbb{N}}$ by setting

$$\mathbf{G}_n = \int_{(n-1)\Delta}^{n\Delta} V_{s-}^{1/2} \, \mathrm{d}L_s. \tag{5.1}$$

Moreover, we shall throughout most of this section presume the following.

***Assumption 5.1.*** *$Y$ (or equivalently $V$) is stationary.*

**Proposition 5.1.** *If Assumption 5.1 holds, then $\mathbf{G}$ is stationary.*

**Proof.** Employing Theorem 4.4 and the same arguments as for [26], Corollary 3.1, shows that $G$ has stationary increments. □

In order to be able to obtain explicit expressions for the moments of $\mathbf{G}$, we need to strengthen Assumption 4.2 as follows.

***Assumption 5.2.*** *Assumption 4.2 is satisfied and, moreover,*

$$E(L_1) = 0 \quad \textit{and} \quad \operatorname{var}(L_1) = (\sigma_W + \sigma_L) I_d \qquad \textit{with } \sigma_W \geq 0.$$

This assumption means that in addition to Assumption 4.2, the Brownian part of $L$ is a scalar multiple of $d$-dimensional standard Brownian motion.

We start by giving conditions for the finiteness of the second moments of $G$, and thus of $\mathbf{G}$, without requiring stationarity and explicit expressions for the moments in the stationary case.

**Proposition 5.2.** *Assume that $E(L_1) = 0$, $E(\|L_1\|^2) < \infty$ and $E(\|Y_0\|) < \infty$. Then $E(\|G_t\|^2) < \infty$ for all $t \in \mathbb{R}^+$.*

*If Assumptions 5.1 and 5.2 are also satisfied, then the stationary sequence $\mathbf{G}$ has the following second-order structure:*

$$E(\mathbf{G}_1) = 0, \tag{5.2}$$

$$\operatorname{var}(\mathbf{G}_1) = (\sigma_L + \sigma_W)\Delta E(V_0), \tag{5.3}$$

$$\operatorname{vec}(\operatorname{var}(\mathbf{G}_1)) = (\sigma_L + \sigma_W)\Delta \mathscr{B}^{-1}(B \otimes I_d + I_d \otimes B)\operatorname{vec}(C), \tag{5.4}$$

$$\operatorname{acov}_{\mathbf{G}}(h) = 0 \qquad \textit{for all } h \in \mathbb{Z}\backslash\{0\}. \tag{5.5}$$

***Remark 5.3.*** (i) In the stationary case, this shows that $\mathbf{G}$ is a white noise sequence.

(ii) Straightforward extensions of the arguments in the proof show that if $Y$ is not stationary, but only (asymptotically) second-order stationary, then $\mathbf{G}$ is (asymptotically) second-order stationary.

For the squared returns $\mathbf{G}\mathbf{G}^* = (\mathbf{G}_n\mathbf{G}_n^*)_{n\in\mathbb{N}}$, we get the following.

**Proposition 5.4.** *Assume that $E(L_1) = 0$, $E(\|L_1\|^4) < \infty$ and $E(\|Y_0\|^2)$ is finite. Then $E(\|G_t\|^4) < \infty$ and, likewise, $E(\|G_tG_t^*\|^2) < \infty$ for all $t \in \mathbb{R}^+$.*

*If Assumptions 5.1 and 5.2 are also satisfied, then the stationary sequence $\mathbf{G}\mathbf{G}^*$ has the following second-order structure:*

$$E(\mathbf{G}_1\mathbf{G}_1^*) = (\sigma_L + \sigma_W)\Delta E(V_0), \tag{5.6}$$

$$\operatorname{acov}_{\mathbf{G}\mathbf{G}^*}(h) = \mathrm{e}^{\mathscr{B}\Delta h}\mathscr{B}^{-1}(I_{d^2} - \mathrm{e}^{-\mathscr{B}\Delta})(\sigma_L + \sigma_W)\operatorname{cov}(\operatorname{vec}(Y_\Delta), \operatorname{vec}(\mathbf{G}_1\mathbf{G}_1^*)) \tag{5.7}$$

*for $h \in \mathbb{N}$.*

Thus, the squared returns $\mathbf{G}\mathbf{G}^*$ have, like an ARMA$(1,1)$ process, a matrix exponentially decreasing autocovariance function from lag one onwards. That such an autocovariance structure is reasonable for financial data can be seen from [34], for instance. In financial data, (quasi-) long-range dependence is frequently encountered, which is often (see, for example, [4]) well modelled by specifying the autocovariance function of the squared increments as the sum of fast and very slowly decaying exponential functions. Since we have a matrix exponential decay, we obtain such a behavior componentwise by appropriate choices for our parameters, with the different rates of the exponential decay being determined by the eigenvalues of $\mathscr{B}$. Additionally we can cover a sinusoidal component.

In the univariate case, [22] obtained, under additional assumptions on $L$, explicit expressions for $\operatorname{var}(\operatorname{vec}(\mathbf{G}_1\mathbf{G}_1^*))$ and $\operatorname{cov}(\operatorname{vec}(Y_\Delta), \operatorname{vec}(\mathbf{G}_1\mathbf{G}_1^*))$. As these are, however, already rather lengthy and complicated formulae, we refrain from calculating these values in our multivariate model.

# 6. Proofs and auxiliary results

In this section, we provide the proofs of our results, along with necessary additional technical results.

## 6.1. Proofs for Section 3

We begin with some matrix analytic results analyzing the Lipschitz properties of the map $V \mapsto V^{1/2} \otimes V^{1/2}$ used in the definition of the MUCOGARCH$(1,1)$ volatility process. We denote by $\|\cdot\|_2$ the operator norm associated with the usual Euclidean norm on $\mathbb{R}^d$.

**Lemma 6.1 ([7], Problem I.6.11).** *For all $A, B \in M_d(\mathbb{R})$, we have*

$$\|A \otimes A - B \otimes B\|_2 \leq 2\max\{\|A\|_2, \|B\|_2\}\|A - B\|_2.$$

*In particular, the mapping $\otimes : M_d(\mathbb{R}) \to M_{d^2}(\mathbb{R}), X \mapsto X \otimes X$ is uniformly Lipschitz on any set of the form $\{x \in M_d(\mathbb{R}) : \|x\| \leq c\}$ with $c > 0$.*

The proof is obvious from the ideas outlined in [7].

**Lemma 6.2 ([7], page 305).** *Let $A, B \in \mathbb{S}_d^+$ and $a > 0$ such that $A, B \geq aI_d$. Then*

$$\|A^{1/2} - B^{1/2}\|_2 \leq \frac{1}{2\sqrt{a}}\|A - B\|_2.$$

*Hence, the mapping $\mathbb{S}_d^+ \to \mathbb{S}_d^+, X \mapsto X^{1/2}$ is uniformly Lipschitz on any set of the form $\{x \in \mathbb{S}_d^+ : x \geq cI_d\} \subset \mathbb{S}_d^{++}$ with $c > 0$.*

For a variant of the above statement see [24], page 557.

**Lemma 6.3.** *Consider the map $F : \mathbb{S}_d^+ \to \mathbb{S}_d^+, X \mapsto X^{1/2} \otimes X^{1/2} = (X \otimes X)^{1/2}$. $F$ is continuous and uniformly Lipschitz on any set of the form $\{x \in \mathbb{S}_d^+ : x \geq cI, \|x\| \leq \tilde{c}\}$ with $c, \tilde{c} > 0$. Moreover, we have that $\|A^{1/2} \otimes A^{1/2}\|_2 = \|A\|_2$ for all $A \in \mathbb{S}_d$.*

**Proof.** The identity $X^{1/2} \otimes X^{1/2} = (X \otimes X)^{1/2}$ is an immediate consequence of basic properties of the tensor product (see [24], Chapter 4) and the continuity of $F$ follows from the continuity of the tensor product and the positive definite square root (see [24], Theorem 6.2.37). The Lipschitz property follows from a combination of the previous two lemmas. Finally, $\|A^{1/2} \otimes A^{1/2}\|_2 = \|A\|_2$ is established by noting that $\|A^{1/2} \otimes A^{1/2}\|_2 = \|A^{1/2}\|_2^2$ (see [7], page 15) and $\|A^{1/2}\|_2 = \|A\|_2^{1/2}$. The latter follows immediately from the fact that $\|A\|_2 = \sqrt{\rho(A^*A)} = \rho(A)$ for all $A \in \mathbb{S}_d^+$. □

Finally, we show that the global Lipschitz property is not satisfied for this map, not even if we restrict it to sets being bounded away from zero.

**Lemma 6.4.** *For the map $F$ defined in the previous lemma, there exists no finite $K \in \mathbb{R}^+$ such that*

$$\|F(x) - F(y)\| \leq K\|x - y\| \tag{6.1}$$

*for all $x, y \in \mathbb{S}_d^{++}$. The same holds for all $x, y \in \{z \in \mathbb{S}_d : z \geq C\}$ with arbitrary $C \in \mathbb{S}_d^{++}$.*

**Proof.** From the following proof, it is clear that we can take $d = 2$ without loss of generality. Let $x = \mathrm{diag}(x_1, x_2)$ and $y = \mathrm{diag}(y_1, y_2)$, with $x_1, x_2, y_1, y_2 \in \mathbb{R}^+\backslash\{0\}$ and $\mathrm{diag}(x_1, x_2)$ being, as usual, the diagonal matrix with diagonal entries $x_1$ and $x_2$. We have $F(x) = \mathrm{diag}(x_1, \sqrt{x_1x_2}, \sqrt{x_1x_2}, x_2)$. Assume that (6.1) is true with a finite $K \in \mathbb{R}^+$.

There is then a finite $k \in \mathbb{R}^+$ such that $|\sqrt{x_1 x_2} - \sqrt{y_1 y_2}| \leq k(|x_1 - y_1| + |x_2 - y_2|)$ for all $x_1, x_2, y_1, y_2 \in \mathbb{R}^+ \backslash \{0\}$. Choosing $x_2 = y_2 = 1$, this gives $|\sqrt{x_1} - \sqrt{y_1}| \leq k|x_1 - y_1|$ for all $x_1, y_1 \in \mathbb{R}^+ \backslash \{0\}$, which is a contradiction to the well-known fact that the square root is not globally Lipschitz on $\mathbb{R}^+ \backslash \{0\}$. Regarding the case $x, y \in \{z \in \mathbb{S}_d : z \geq C\}$, we can, without loss of generality, restrict ourselves to $C = cI_d$ with $c \in \mathbb{R}^+ \backslash \{0\}$. Choosing $x_2 = y_2$, $x_1 = 9c$ and $y_1 = 4c$ gives $|\sqrt{cx_2}| \leq 5kc$. As $x_2$ can be taken arbitrarily large, this is a contradiction. □

In the following, we use the fact that any stochastic differential equation defined on an open set which has locally Lipschitz coefficients growing at most linearly has a unique solution until the first time the open set is left or its boundary is reached. This result follows along the same lines as the usual existence results for SDEs with locally Lipschitz coefficients defined on $\mathbb{R}^d$ (see, for example, [33] or [36], Theorem V.38). Alternatively, a proof for open sets of the type relevant below can be found in [40], Section 6.7, which uses only the standard existence and uniqueness results for SDEs on $\mathbb{R}^d$ with globally Lipschitz coefficients and orthogonal projections.

**Proof of Theorem 3.2.** Define the maps $F$ and $G$ by $F(\mathrm{vec}(y)) = (I_d \otimes B + B \otimes I_d)\,\mathrm{vec}(y)$ and $G(y) = (A \otimes A)((C + y)^{1/2} \otimes (C + y)^{1/2})$. The SDE (3.4) can then be written as

$$\mathrm{d}\,\mathrm{vec}(Y_t) = F(\mathrm{vec}(Y_{t-}))\,\mathrm{d}t + G(Y_{t-})\,\mathrm{d}\,\mathrm{vec}([L, L]^{\mathfrak{d}}_t). \tag{6.2}$$

Moreover, we define the set $U_{C,\varepsilon} = \{x \in \mathbb{S}_d : x > -\varepsilon I_d\}$ for some $\varepsilon$ with $0 < \varepsilon < \min \sigma(C)$. Then the set $U_{C,\varepsilon}$ (and thus $\mathrm{vec}(U_{C,\varepsilon})$) is open and for each $x \in U_{C,\varepsilon}$, we have $x + C > (\min \sigma(C) - \varepsilon) I_d \in \mathbb{S}_d^{++}$. Since the foregoing results imply that $G$ is locally Lipschitz on $U_{C,\varepsilon}$ and has linear growth (a function $f$ has linear growth if $\|f(x)\|^2 \leq C(1 + \|x\|^2)$), standard results on the existence of solutions of SDEs give that (6.2) has a unique locally bounded solution $(Y_t)_{t \in \mathbb{R}^+}$ with initial value $Y_0$, provided it can be ensured that every solution does not leave the set $U_{C,\varepsilon}$ or touch its boundary. However, it is easy to see that every solution must satisfy $Y_t \geq \mathrm{e}^{Bt} Y_0 \mathrm{e}^{B^* t}$ since all jumps $A V_{t-}^{1/2} \Delta[L, L]^{\mathfrak{d}}_t V_{t-}^{1/2} A^*$ are positive semidefinite and between jumps, $Y$ follows the deterministic differential equation $\mathrm{d}Y_t = (BY_{t-} + Y_{t-} B^*)\,\mathrm{d}t$ uniquely solved by $Y_t = \mathrm{e}^{Bt} Y_0 \mathrm{e}^{B^* t}$, so any solution necessarily stays in $\mathbb{S}_d^+$.

The finite variation property is clear since time $t$ and $[L, L]^{\mathfrak{d}}$ are of finite variation. □

**Proof of Theorem 3.6.** Define $M_t = \int_0^t A(C + Y_{s-})^{1/2}\,\mathrm{d}[L, L]^{\mathfrak{d}}_s\,(C + Y_{s-})^{1/2} A^*$. Then $M$ is $\mathbb{S}_d^+$-increasing and of finite variation and $Y$ obviously solves the stochastic differential equation $\mathrm{d}X_t = (BX_{t-} + X_{t-} B^*)\,\mathrm{d}t + \mathrm{d}M_t (*)$. Standard theory implies that this differential equation has a unique solution, and the same elementary calculations as for Ornstein–Uhlenbeck processes show that the solution of $(*)$ with initial value $Y_0$, which is necessarily equal to $Y$, is given by

$$\mathrm{e}^{Bt} Y_0 \mathrm{e}^{B^* t} + \int_0^t \mathrm{e}^{B(t-s)} A\,\mathrm{d}M_s A^* \mathrm{e}^{B^*(t-s)}$$

$$= \mathrm{e}^{Bt} Y_0 \mathrm{e}^{B^* t} + \int_0^t \mathrm{e}^{B(t-s)} A (C + Y_{s-})^{1/2} \, \mathrm{d}[L, L]^{\mathfrak{d}}_s \, (C + Y_{s-})^{1/2} A^* \mathrm{e}^{B^*(t-s)}. \qquad \square$$

## 6.2. Proofs for Section 4

### *6.2.1. Proofs for Section 4.1*

The univariate COGARCH(1, 1) bounds will first be shown for processes driven by compound Poisson processes and then for the general case using an approximation by compound Poisson processes which is of interest in its own right, as it provides, for instance, a possible approximation scheme to be used in simulations.

To see that the processes defined in the following are indeed univariate COGARCH processes, we need the following general lemma.

**Lemma 6.5.** *Let $(L_t)_{t \in \mathbb{R}^+}$ be a driftless Lévy subordinator. There then exists a Lévy process $(\overline{L}_t)_{t \in \mathbb{R}^+}$ in $\mathbb{R}$ such that $L_t = [\overline{L}, \overline{L}]^{\mathfrak{d}}_t$ for all $t \in \mathbb{R}^+$.*

**Proof.** Denote the jump measure associated with $L$ by $\mu_L$, that is, $L_t = \int_0^t \int_{\mathbb{R}^+} x \mu_L(\mathrm{d}s, \mathrm{d}x)$, and denote its Lévy measure by $\nu_L$. Let

$$\overline{L}_t = \int_0^t \int_{0 < x \le 1} \sqrt{x} (\mu_L(\mathrm{d}s, \mathrm{d}x) - \mathrm{d}s \, \nu_L(\mathrm{d}x)) + \int_0^t \int_{x > 1} \sqrt{x} \mu_L(\mathrm{d}s, \mathrm{d}x),$$

noting that the existence of the first integral follows from the finiteness of $\int_{0 < x \le 1} \sqrt{x}^2 \nu_L(\mathrm{d}x) = \int_{0 < x \le 1} x \nu_L(\mathrm{d}x)$ since $L$ is of finite variation. Now $(\overline{L}_t)_{t \in \mathbb{R}^+}$ is a Lévy process and $L_t = [\overline{L}, \overline{L}]^{\mathfrak{d}}_t$. $\square$

The following elementary results regarding the norm $\| \cdot \|_{B,S}$ are straightforward to obtain.

**Lemma 6.6.** *It holds that $\|S \otimes S\|_{B,S} = \|S\|_2^2$ and $\|S^{-1} \otimes S^{-1}\|_{B,S} = \|S^{-1}\|_2^2$. Moreover,*

$$\|x\|_{B,S} \le \|S^{-1}\|_2^2 \|x\|_2 \quad \text{and} \quad \|x\|_2 \le \|S\|_2^2 \|x\|_{B,S} \qquad \text{for all } x \in \mathbb{R}^{d^2},$$
$$\|X\|_{B,S} \le \|S\|_2^2 \|S^{-1}\|_2^2 \|X\|_2 \quad \text{and} \quad \|X\|_2 \le \|S\|_2^2 \|S^{-1}\|_2^2 \|X\|_{B,S} \qquad \text{for all } X \in M_{d^2}(\mathbb{R}).$$

Note that a very similar norm has also been used in [10].

We can now analyze the norms of compound Poisson driven MUCOGARCH volatility processes. Recall that in the univariate case, the MUCOGARCH volatility process $Y$ is just a deterministically scaled version of the COGARCH volatility process $\mathbf{Y}$ defined in [10].

**Proof of Theorem 4.1 if $L$ is compound Poisson.** Using Lemma 6.5, it is clear that the process $(y_t)_{t \in \mathbb{R}}$ defined in Theorem 4.1 is a univariate MUCOGARCH(1, 1) process.

Let $\Gamma_1$ be the time of the first jump of $L$ and let $t\in[0,\Gamma_1)$. Since $\|\mathrm{e}^{(I_d\otimes B+B\otimes I_d)t}\|_{B,S}=\mathrm{e}^{2\lambda t}$, it holds that

$$\begin{aligned}\|\operatorname{vec}(Y_t)\|_{B,S}&=\|\mathrm{e}^{(I_d\otimes B+B\otimes I_d)t}\operatorname{vec}(Y_0)\|_{B,S}\le\|\mathrm{e}^{(I_d\otimes B+B\otimes I_d)t}\|_{B,S}\|\operatorname{vec}(Y_0)\|_{B,S}\\&=\mathrm{e}^{2\lambda t}y_0=y_t.\end{aligned}$$

Thus, (4.2) is shown for all $t\in[0,\Gamma_1)$. At time $\Gamma_1$, we have

$$\begin{aligned}&\|\operatorname{vec}(Y_{\Gamma_1})\|_{B,S}\\&\quad=\|\operatorname{vec}(Y_{\Gamma_1-})+(A\otimes A)((C+Y_{\Gamma_1-})^{1/2}\otimes(C+Y_{\Gamma_1-})^{1/2})\operatorname{vec}(\Delta L_{\Gamma_1}(\Delta L_{\Gamma_1})^*)\|_{B,S}\\&\quad\le y_{\Gamma_1-}+\|A\otimes A\|_{B,S}\|(C+Y_{\Gamma_1-})^{1/2}\otimes(C+Y_{\Gamma_1-})^{1/2}\|_{B,S}\|\operatorname{vec}(\Delta L_{\Gamma_1}(\Delta L_{\Gamma_1})^*)\|_{B,S}\\&\quad\le y_{\Gamma_1-}+\|A\otimes A\|_{B,S}\|S\|_2^2\|S^{-1}\|_2^2\|(C+Y_{\Gamma_1-})^{1/2}\otimes(C+Y_{\Gamma_1-})^{1/2}\|_2\Delta\tilde{L}_{\Gamma_1}\\&\quad\le y_{\Gamma_1-}+\|A\otimes A\|_{B,S}\|S\|_2^2\|S^{-1}\|_2^2(\|C\|_2+\|Y_{\Gamma_1-}\|_2)\Delta\tilde{L}_{\Gamma_1}\\&\quad\le y_{\Gamma_1-}+\|A\otimes A\|_{B,S}\|S\|_2^2\|S^{-1}\|_2^2K_{2,B}(K_{2,B}^{-1}\|C\|_2+\|\operatorname{vec}(Y_{\Gamma_1-})\|_{B,S})\Delta\tilde{L}_{\Gamma_1}=y_{\Gamma_1},\end{aligned}$$

which establishes (4.2) for $t=\Gamma_1$. Iterating these arguments shows (4.2) for all $t\in\mathbb{R}^+$.

The first inequality for $K_{2,B}$ follows from Lemma 6.6 and the second one by [23], page 314. □

In order to extend Theorem 4.1 to MUCOGARCH processes driven by general Lévy processes, we need to show that we can approximate a MUCOGARCH volatility process by approximating the driving Lévy process. The following result is very similar to [10], Lemma 8.2. However, we need to give a detailed proof since the standard results cannot be applied due to the fact that we have only locally Lipschitz coefficients.

**Proposition 6.7.** *Let $Y$ be a* MUCOGARCH *volatility process with $C\in\mathbb{S}_d^{++}$ and $Y_0\in\mathbb{S}_d^+$, driven by a Lévy process $L$ in $\mathbb{R}^d$, and let $(\varepsilon_n)_{n\in\mathbb{N}}$ be a monotonically decreasing sequence in $\mathbb{R}^+\backslash\{0\}$ with $\lim_{n\to\infty}\varepsilon_n=0$. Define compound Poisson Lévy processes $L_n$ by $L_{n,t}=\int_0^t\int_{\mathbb{R}^d,\|x\|\ge\varepsilon_n}x\mu_L(\mathrm{d}s,\mathrm{d}x)$ for $n\in\mathbb{N}$ and associated* MUCOGARCH *volatility processes $Y_n$ by*

$$\begin{aligned}\mathrm{d}Y_{n,t}&=(BY_{n,t-}+Y_{n,t-}B^*)\,\mathrm{d}t+A(C+Y_{n,t-})^{1/2}\,\mathrm{d}[L_n,L_n]_t^{\mathfrak{d}}\,(C+Y_{n,t-})^{1/2}A^*,\\Y_{n,0}&=Y_0.\end{aligned}\tag{6.3}$$

*Then $Y_n\to Y$ as $n\to\infty$ almost surely uniformly on compacts.*

**Proof.** First, observe that $[L_n,L_n]_t^{\mathfrak{d}}=\int_0^t\int_{\|x\|\ge\varepsilon_n}xx^*\mu_L(\mathrm{d}s,\mathrm{d}x)$ implies that $[L_n,L_n]^{\mathfrak{d}}\to[L,L]^{\mathfrak{d}}$, as $n\to\infty$ a.s. uniformly on compacts, and that $[L_n,L_n]^{\mathfrak{d}}$ is monotonically increasing in $n$.

Since all processes involved are of finite variation, we can prove the claim with a pathwise approach. So, fix $\omega \in \Omega$ and thereby one path. Let $T \in \mathbb{R}^+$ be arbitrary. The Gronwall inequality (see [36], Exercise 15, page 358) shows that

$$\|Y_{n,t}\|_2 \leq \left(\|Y_0\|_2 + \|A\|_2^2\|C\|_2 \int_0^T \int_{\mathbb{R}^d} \|x\|_2^2 \mu_{L_n}(\mathrm{d}s, \mathrm{d}x)\right) \mathrm{e}^{\|A\|_2^2 \int_0^t \int_{\mathbb{R}^d} \|x\|_2^2 \mu_{L_n}(\mathrm{d}s,\mathrm{d}x) + 2\|B\|_2 t}$$

$$\leq \left(\|Y_0\|_2 + \|A\|_2^2\|C\|_2 \int_0^T \int_{\mathbb{R}^d} \|x\|_2^2 \mu_{L}(\mathrm{d}s, \mathrm{d}x)\right) \mathrm{e}^{\|A\|_2^2 \int_0^T \int_{\mathbb{R}^d} \|x\|_2^2 \mu_{L}(\mathrm{d}s,\mathrm{d}x) + 2\|B\|_2 T},$$

$$\|Y_t\|_2 \leq \left(\|Y_0\|_2 + \|A\|_2^2\|C\|_2 \int_0^T \int_{\mathbb{R}^d} \|x\|_2^2 \mu_{L}(\mathrm{d}s, \mathrm{d}x)\right) \mathrm{e}^{\|A\|_2^2 \int_0^T \int_{\mathbb{R}^d} \|x\|_2^2 \mu_{L}(\mathrm{d}s,\mathrm{d}x) + 2\|B\|_2 T}$$

for all $t \in [0,T]$. Since $Y_t \geq \mathrm{e}^{Bt} Y_0 \mathrm{e}^{B^* t}$, $Y_{n,t} \geq \mathrm{e}^{Bt} Y_0 \mathrm{e}^{B^* t}$ and $Y_0$ is positive semidefinite, $C + Y$ and $(C + Y_n)_{n\in\mathbb{N}}$ all remain in one common compact set in $\mathbb{S}_d^{++}$ on $[0,T]$. Thus, when considering (3.4) and (6.3), we can regard the coefficients of these SDEs as being globally Lipschitz with a common Lipschitz coefficient. Thus, [36], Corollary, page 261 after Theorem v.11, implies that $Y_n(\omega) \to Y(\omega)$ uniformly on $[0,T]$. Note that, formally, the result of [36] is applied on the probability space given by the set $\{\omega\}$, the trivial $\sigma$-algebra $\{\{\omega\}, \varnothing\}$ (which also gives the filtration) and the Dirac measure with respect to $\omega$.

Since $\omega \in \Omega$ and $T \in \mathbb{R}^+$ were arbitrary, this completes the proof. □

**Proof of Theorem 4.1 for general *L*.** Let $(Y_n)_{n\in\mathbb{N}}$ be the sequence of compound Poisson driven MUCOGARCH(1,1) processes converging a.s. on compacts to $Y$ constructed in the last proposition. For $n \in \mathbb{N}$, denote by $y_n$ the univariate MUCOGARCH(1,1) processes with $\|\operatorname{vec}(Y_{n,t})\|_{B,S} \leq y_{n,t}$ for all $t \in \mathbb{R}^+$. Then $y_{n,t} + K_{2,B}^{-1}\|C\|$ is a univariate COGARCH(1,1) volatility process as defined in [26], where it is denoted by $\sigma_{t+}^2$. Since we only add more jumps in $[L_n, L_n]^{\mathfrak{d}}$ when we increase $n$, it is straightforward to see from equations (3.3) and (3.4) in [26] that $y_{n+l,t} \geq y_{n,t}$ for all $n, l \in \mathbb{N}$ and $t \in \mathbb{R}^+$. Moreover, defining the process $y$ by

$$\begin{aligned} \mathrm{d}y_t &= 2\lambda y_{t-}\,\mathrm{d}t + \|S\|_2^2\|S^{-1}\|_2^2 K_{2,B}\|A \otimes A\|_{B,S}\left(\frac{\|C\|_2}{K_{2,B}} + y_{t-}\right)\mathrm{d}\tilde{L}_t, \\ y_0 &= \|\operatorname{vec}(Y_0)\|_{B,S}, \end{aligned} \tag{6.4}$$

with $\tilde{L}_t := \int_0^t \int_{\mathbb{R}^d} \|\operatorname{vec}(xx^*)\|_{B,S}\mu_L(\mathrm{d}s, \mathrm{d}x)$, the same argument implies that $y_{n,t} \leq y_t$ for all $n \in \mathbb{N}$ and $t \in \mathbb{R}^+$. Note that $(\tilde{L}_t)_{t\in\mathbb{R}^+}$ is a well-defined Lévy process, because there is a $K > 0$ such that

$$\int_{\mathbb{R}^d} (\|\operatorname{vec}(xx^*)\|_{B,S} \wedge 1)\nu_L(\mathrm{d}x) \leq K \int_{\mathbb{R}^d} (\|xx^*\|_2 \wedge 1)\nu_L(\mathrm{d}x) = K \int_{\mathbb{R}^d} (\|x\|_2^2 \wedge 1)\nu_L(\mathrm{d}x) < \infty.$$

Passing to the limit $n \to \infty$ in $\|\operatorname{vec}(Y_{n,t})\|_{B,S} \leq y_{n,t} \leq y_t$ establishes $\|\operatorname{vec}(Y_t)\|_{B,S} \leq y_t$ for all $t \in \mathbb{R}^+$. □

**Proof of Proposition 4.3.** Let $y$ be the process constructed in Theorem 4.1. It then suffices to show that $E(y_t^k) < \infty$ and that this is locally bounded in $t$. By construction, $E(y_0^k)$ is finite. Moreover, let $\bar{L}$ be the Lévy process constructed in Lemma 6.5 such that $\tilde{L}_t = [\bar{L}, \bar{L}]_t^{\mathfrak{d}}$. The finiteness of $E(\|L_1\|^{2k})$ implies that $\int_{\mathbb{R}^d} \|x\|_2^{2k} \nu_L(\mathrm{d}x) = \int_{\mathbb{R}^d} \|xx^*\|_2^k \nu_L(\mathrm{d}x) < \infty$ (with $\|\cdot\|_2$ denoting the Euclidean norm in the first integral and the associated operator norm in the second integral). Since the finiteness of the integrals is independent of the particular norm used, it follows that $\int_{\mathbb{R}^d} \|\operatorname{vec}(xx^*)\|_{B,S}^k \nu_L(\mathrm{d}x) = \int_{\mathbb{R}} |x|^k \nu_{\tilde{L}}(\mathrm{d}x) = \int_{\mathbb{R}} |x|^{2k} \nu_{\bar{L}}(\mathrm{d}x) < \infty$. Hence, $E(|\bar{L}_1|^{2k})$ is finite and using the results of [26], Section 4, as in the proof of [10], Proposition 4.1, completes the proof of this proposition. □

### *6.2.2. Proofs for Section 4.2*

In order to show the existence of a stationary distribution of the MUCOGARCH volatility process $Y$, we need to recall a result from the theory of weak convergence. For more details and the relevant background, we refer to any of the standard texts (for example, [8, 25]). Below, we denote by $\mathcal{M}_1(E)$ the set of all probability measures on the Borel $\sigma$-algebra of a Polish space $E$.

The following theorem on the existence of a stationary distribution for a Markov process is referred to as the "Krylov–Bogoliubov existence theorem" in the literature. For a proof, see [13], Section 3.1, or [37], Theorem 4.6.

**Theorem 6.8.** *Let $E$ be a Polish space and $(P_s)_{s\in\mathbb{R}^+}$ the transition semigroup of an $E$-valued weak Feller Markov process. Assume that there is an $\eta \in \mathcal{M}_1(E)$ such that the set $\{P_t^*\eta : t \in \mathbb{R}^+\}$ is tight. There then exists a $\mu \in \mathcal{M}_1(E)$ such that $P_t^*\mu = \mu$ for all $t \in \mathbb{R}^+$, that is, $\mu$ is an invariant measure for $(P_s)_{s\in\mathbb{R}^+}$ or a stationary distribution for the Markov process, respectively, and $\mu$ is in the closed (with respect to weak convergence) convex hull of $\{P_t^*\eta : t \in \mathbb{R}^+\}$.*

Above, $P_t^* : \mathcal{M}_1(E) \to \mathcal{M}_1(E)$ denotes the operator given by $P_t^*\mu(U) = \int_E P_s(x, U)\mu(\mathrm{d}x)$ for any Borel set $U$ where $P_t(x, U)$ is the transition probability of the Markov process from the initial state $x$ to the set $U$ at time $t \in \mathbb{R}^+$.

**Proof.** Proof of Theorem 4.5 Let $\lambda, \tilde{L}$ be defined as in Theorem 4.1 and $\bar{L}$ be the Lévy process constructed in Lemma 6.5 such that $\tilde{L}_t = [\bar{L}, \bar{L}]_t^{\mathfrak{d}}$. Then

$$\int_{\mathbb{R}^d} \log(1 + \alpha_1 \|\operatorname{vec}(yy^*)\|_{B,S}) \nu_L(\mathrm{d}y) = \int_{\mathbb{R}^d} \log(1 + \alpha_1 y^2) \nu_{\bar{L}}(\mathrm{d}y)$$

and thus [10], Theorem 3.1, (see also [26], Theorem 3.1), shows that the process $y$ satisfying (4.1) converges in distribution to a distribution concentrated on $\mathbb{R}^+$. Assume now that $y_0$ has this stationary probability distribution and is independent of $(L_s)_{s\in\mathbb{R}^+}$. Setting $Y_0 = \frac{y_0}{\|\operatorname{vec}(I_d)\|_{B,S}} I_d$ gives an initial value for the MUCOGARCH volatility process that is independent of $L$ and, moreover, $\|\operatorname{vec}(Y_0)\|_{B,S} = y_0$. Thus, the process $y$ satisfying $\|\operatorname{vec}(Y_t)\|_{B,S} \le y_t$ for all $t \in \mathbb{R}^+$ (see Theorem 4.1) is stationary. Since for every $K>0$, the set $\{x \in \mathbb{S}_d : \|x\| \le K\}$ is compact in $\mathbb{S}_d^+$, $P(\|Y_t\|_{B,S} \le K) \ge P(y_t \le K)$ and $y$ is stationary with a stationary distribution concentrated on $\mathbb{R}^+$, it follows that the set $\{\mathscr{L}(Y_t) : t \in \mathbb{R}^+\}$ of laws $\mathscr{L}(Y_t)$ of $Y_t$ forms a tight subset of $\mathcal{M}_1(\mathbb{S}_d^+)$. Therefore, Theorem 6.8 combined with Theorem 4.4 implies that there exists a stationary distribution $\mu \in \mathcal{M}_1(\mathbb{S}_d^+)$ for the MUCOGARCH volatility process $Y$ such that $\mu$ is in the closed convex hull of $\{\mathscr{L}(Y_t) : t \in \mathbb{R}^+\}$.

If (4.4) holds for some $k \in \mathbb{N}$, [10], Proposition 4.1 (see also [26], Section 4), shows that the stationary distribution of $y$ has a finite $k$th moment. This, in turn, implies that $E(\|Y_t\|^k) \le c$ for some finite $c \in \mathbb{R}^+$ and all $t \in \mathbb{R}^+$. Hence, $\int_{\mathbb{S}_d^+} \|x\|^k \mu(\mathrm{d}x) < \infty$ because $\mu$ is in the closed convex hull of $\{\mathscr{L}(Y_t) : t \in \mathbb{R}^+\}$. □

### *6.2.3. Proofs for Section 4.3*

In the following calculations of moments of $Y$, we often use the fact that the stochastic continuity of $L$ implies $E(Y_{t-}) = E(Y_t)$, $\operatorname{var}(\operatorname{vec}(Y_{t-})) = \operatorname{var}(\operatorname{vec}(Y_t))$ and similar results. Moreover, the following version of the so-called compensation formula is needed.

**Lemma 6.9.** *Assume that $(X_t)_{t\in\mathbb{R}^+}$ is an adapted cadlag $M_d(\mathbb{R})$-valued process satisfying $E(\|X_t\|) < \infty$ for all $t \in \mathbb{R}^+$, $t \mapsto E(\|X_t\|)$ is locally bounded and $(L_t)_{t\in\mathbb{R}^+}$ is a driftless pure jump Lévy process in $\mathbb{R}^d$ of finite variation with finite expectation $E(\|L_1\|)$. Then $E(\int_0^t X_{s-}\,\mathrm{d}L_s) = \int_0^t E(X_{s-})E(L_1)\,\mathrm{d}s$ for $t \in \mathbb{R}^+$.*

**Proof.** Since $L_t = \int_0^t \int_{\mathbb{R}^d} z\mu_L(\mathrm{d}s,\mathrm{d}z)$, the compensation formula (see [28], Section 4.3.2) implies that

$$\begin{aligned} E\biggl(\int_0^t X_{s-}\,\mathrm{d}L_s\biggr) &= E\biggl(\int_0^t \int_{\mathbb{R}^d} X_{s-} z\mu_L(\mathrm{d}s,\mathrm{d}z)\biggr) = E\biggl(\int_0^t X_{s-} \int_{\mathbb{R}^d} z\nu_L(\mathrm{d}z)\,\mathrm{d}s\biggr) \\ &= E\biggl(\int_0^t X_{s-} E(L_1)\,\mathrm{d}s\biggr). \end{aligned}$$

Observing that $\int_0^t E(\|X_{s-}\|)\|E(L_1)\|\,\mathrm{d}s$ is finite for every $t \in \mathbb{R}^+$, an application of Fubini's theorem completes the proof. □

**Proof of Proposition 4.7.** Consider the case $k=1$ first. Elementary arguments give that

$$\begin{aligned} \|Y_t\|_2 \le \|Y_0\|_2 &+ \|A\|_2^2 \|C\|_2 \int_0^t \int_{\mathbb{R}^d} \|x\|_2^2 \mu_L(\mathrm{d}s,\mathrm{d}x) + 2\|B\|_2 \int_0^t \|Y_{s-}\|_2\,\mathrm{d}s \\ &+ \|A\|_2^2 \int_0^t \int_{\mathbb{R}^d} \|Y_{s-}\|_2 \|x\|_2^2 \mu_L(\mathrm{d}s,\mathrm{d}x). \end{aligned}$$

Using stochastic continuity, the compensation formula and the observation that $\int_{\mathbb{R}^d}\|x\|_2^2\nu_L(\mathrm{d}x)<\infty$ due to the finiteness of $E(\|L_1\|_2^2)$, this implies that

$$E(\|Y_t\|_2)\le E(\|Y_0\|_2)+\|A\|_2^2\|C\|_2T\int_{\mathbb{R}^d}\|x\|_2^2\nu_L(\mathrm{d}x) + \left(2\|B\|_2+\|A\|_2^2\int_{\mathbb{R}^d}\|x\|_2^2\nu_L(\mathrm{d}x)\right)\int_0^t E(\|Y_s\|_2)\,\mathrm{d}s$$

for all $t\in[0,T]$ and any $T\in\mathbb{R}^+$. The Gronwall lemma thus shows that $E(\|Y_t\|_2)$ is finite and bounded for $t\in[0,T]$. Since $T$ was arbitrary, this completes the proof for this case.

In the case $k\ge 2$, we obtain from [36], Theorem V.66, and the elementary inequality $|a+b|^k\le 2^{k-1}(|a|^k+|b|^k)$ for all $a,b\in\mathbb{R}$ that there exists a constant $K_k\in\mathbb{R}^+$ such that

$$E(\|Y_t\|_2^k)\le K_k\left(E(\|Y_0\|_2^k)+\|A\|_2^{2k}\|C\|_2^k t+(2^k\|B\|_2^k+\|A\|_2^{2k})\int_0^t E(\|Y_s\|_2^k)\,\mathrm{d}s\right).$$

Using the Gronwall lemma and arguing analogously to the case $k=1$ now completes the proof. □

**Proof of Theorem 4.8.** Proposition 4.7 ensures the finiteness and local boundedness of the first absolute moment needed in the following. From the defining stochastic differential equation (3.4), we have

$$Y_t=Y_0+\int_0^t(BY_{s-}+Y_{s-}B^*)\,\mathrm{d}s+\int_0^t A(Y_{s-}+C)^{1/2}\,\mathrm{d}[L,L]_s^{\mathfrak{d}}(Y_{s-}+C)^{1/2}A^*.$$

Therefore,

$$E(Y_t)=E(Y_0)+\int_0^t(BE(Y_s)+E(Y_s)B^*)\,\mathrm{d}s+\sigma_L\int_0^t AE(Y_s+C)A^*\,\mathrm{d}s, \tag{6.5}$$

using a Fubini argument, stochastic continuity, the variant of the compensation formula given in Lemma 6.9 and the observation that $E([L,L]_1^{\mathfrak{d}})=\int_{\mathbb{R}^d}xx^*\nu(\mathrm{d}x)=\mathrm{var}(L_1^{\mathfrak{d}})$ is implied by equation (2.1). Thus,

$$\begin{aligned}
&\mathrm{vec}\left(E\left(\int_0^t A(Y_{s-}+C)^{1/2}\,\mathrm{d}[L,L]_s^{\mathfrak{d}}\,(Y_{s-}+C)^{1/2}A^*\right)\right)\\
&\quad=E\left(\int_0^t(A\otimes A)((Y_{s-}+C)^{1/2}\otimes(Y_{s-}+C)^{1/2})\,\mathrm{d}\,\mathrm{vec}([L,L]_s^{\mathfrak{d}})\right)\\
&\quad=\int_0^t(A\otimes A)E((Y_{s-}+C)^{1/2}\otimes(Y_{s-}+C)^{1/2})\,\mathrm{vec}(E([L,L]_1^{\mathfrak{d}}))\,\mathrm{d}s\\
&\quad=\sigma_L\int_0^t(A\otimes A)E((Y_{s-}+C)^{1/2}\otimes(Y_{s-}+C)^{1/2}\,\mathrm{vec}(I_d))\,\mathrm{d}s
\end{aligned}$$

$$= \sigma_L \operatorname{vec}\left( \int_0^t AE(Y_{s-} + C)A^* \, \mathrm{d}s \right).$$

Equation (6.5) therefore implies the following differential equation after vectorizing:

$$\frac{\mathrm{d}}{\mathrm{d}t} E(\operatorname{vec}(Y_t)) = \mathscr{B} E(\operatorname{vec}(Y_t)) + \sigma_L (A \otimes A) \operatorname{vec}(C).$$

Solving this ODE establishes (i).

Turning to (ii), the assumed second-order stationarity and (6.5) imply that

$$BE(Y_0) + E(Y_0)B^* + \sigma_L A(E(Y_0) + C)A^* = 0.$$

The rest is just a matter of rewriting this linear equation. □

**Proof of Lemma 4.10.** From Assumption 4.2, we have that

$$\sigma_L \| \operatorname{vec}(I_d) \|_{B,S} = \left\| \int_{\mathbb{R}^d} \operatorname{vec}(xx^*) \nu_L(\mathrm{d}x) \right\|_{B,S} \le \int_{\mathbb{R}^d} \| \operatorname{vec}(xx^*) \|_{B,S} \nu_L(\mathrm{d}x). \tag{6.6}$$

For $k = 1$, condition (4.4) becomes $\|S\|_2^2 \|S^{-1}\|_2^2 K_{2,B} \|A \otimes A\|_{B,S} \int_{\mathbb{R}^d} \| \operatorname{vec}(xx^*) \|_{B,S} \nu_L(\mathrm{d}x) < -2\lambda$. Using (6.6), $\|S\|_2 \|S^{-1}\|_2 \ge 1$ and the fact that $K_{2,B} \| \operatorname{vec}(I_d) \|_{B,S} \ge \|I_d\|_2 = 1$ due to the definition of $K_{2,B}$, one obtains

$$\sigma_L \|A \otimes A\|_{B,S} = \sigma_L \|(S^{-1} \otimes S^{-1})(A \otimes A)(S \otimes S)\|_2 < -2\lambda.$$

Let $\mu$ now be any eigenvalue of $\mathscr{B}$ and note that $(S^{-1} \otimes S^{-1})(B \otimes I_d + I_d \otimes B)(S \otimes S)$ is diagonal. Thus, the Bauer–Fike theorem (see [23], Theorem 6.3.2 and its proof, for instance) gives that there exists a $\tilde{\mu} \in \sigma(B \otimes I_d + I_d \otimes B)$ such that

$$\begin{aligned} |\Re(\mu) - \Re(\tilde{\mu})| &\le |\mu - \tilde{\mu}| \le \|(S^{-1} \otimes S^{-1})(\mathscr{B} - B \otimes I_d - I_d \otimes B)(S \otimes S)\| \\ &= \sigma_L \|(S^{-1} \otimes S^{-1})(A \otimes A)(S \otimes S)\|_2 < -2\lambda. \end{aligned}$$

Hence, $\Re(\mu) < \max\{\Re(\tilde{\mu}) : \tilde{\mu} \in \sigma(B \otimes I_d + I_d \otimes B)\} - 2\lambda = 0$ because the maximum equals $2\lambda$ due to $\sigma(B \otimes I_d + I_d \otimes B) = \sigma(B) + \sigma(B)$ and the definition of $\lambda = \max(\Re(\sigma(B)))$. Therefore, $\sigma(\mathscr{B}) \subset (-\infty, 0) + \mathrm{i}\mathbb{R}$ and $\mathscr{B}$ is invertible. □

**Proof of Theorem 4.11.** We only prove (i), because the proof of (ii) proceeds along the same lines, noting that the finiteness is ensured by Proposition 4.7.

The equality $\operatorname{acov}_Y(\cdot) = \operatorname{acov}_V(\cdot)$ is obvious. Due to the second-order stationarity, we have

$$\begin{aligned} \operatorname{acov}_Y(h) = \operatorname{cov}\Bigg( \operatorname{vec}\bigg( Y_0 &+ \int_0^h (BY_{s-} + Y_{s-}B^*) \, \mathrm{d}s \\ &+ \int_0^h A(Y_{s-} + C)^{1/2} \, \mathrm{d}[L, L]_s^{\mathfrak{d}} \, (Y_{s-} + C)^{1/2} A^* \bigg), \operatorname{vec}(Y_0) \Bigg) \end{aligned}$$

$$
\begin{aligned}
&= \mathrm{var}(\mathrm{vec}(Y_0)) + E\left(\int_0^h (B\otimes I_d + I_d\otimes B)\,\mathrm{vec}(Y_{s-})\,\mathrm{vec}(Y_0)^*\,\mathrm{d}s\right)\\
&\quad - E\left(\int_0^t (B\otimes I_d + I_d\otimes B)\,\mathrm{vec}(Y_{s-})\,\mathrm{d}s\right)E(\mathrm{vec}(Y_0)^*)\\
&\quad + E\left(\int_0^h (A\otimes A)((Y_{s-}+C)^{1/2}\otimes(Y_{s-}+C)^{1/2})\,\mathrm{d}\,\mathrm{vec}([L,L]_s^{\mathfrak{d}})\ \mathrm{vec}(Y_0)^*\right)\\
&\quad - E\left(\int_0^h (A\otimes A)((Y_{s-}+C)^{1/2}\otimes(Y_{s-}+C)^{1/2})\,\mathrm{d}\,\mathrm{vec}([L,L]_s^{\mathfrak{d}})\right)E(\mathrm{vec}(Y_0)^*)\\
&= \mathrm{var}(\mathrm{vec}(Y_0)) + \int_0^h (B\otimes I_d + I_d\otimes B)E(\mathrm{vec}(Y_s)\,\mathrm{vec}(Y_0)^*)\,\mathrm{d}s\\
&\quad - \int_0^t (B\otimes I_d + I_d\otimes B)E(\mathrm{vec}(Y_s))\,\mathrm{d}s\,E(\mathrm{vec}(Y_0)^*)\\
&\quad + \sigma_L\int_0^h (A\otimes A)E(((Y_{s-}+C)^{1/2}\otimes(Y_{s-}+C)^{1/2})\,\mathrm{vec}(I_d)\,\mathrm{vec}(Y_0)^*)\,\mathrm{d}s\\
&\quad - \sigma_L\int_0^h (A\otimes A)E(((Y_{s-}+C)^{1/2}\otimes(Y_{s-}+C)^{1/2})\,\mathrm{vec}(I_d))\,\mathrm{d}sE(\mathrm{vec}(Y_0)^*)\\
&= \sigma_L\int_0^h (A\otimes A)E(\mathrm{vec}(Y_s+C)\,\mathrm{vec}(Y_0)^*)\,\mathrm{d}s\\
&\quad - \sigma_L\int_0^h (A\otimes A)E(\mathrm{vec}(Y_s+C))E(\mathrm{vec}(Y_0)^*)\,\mathrm{d}s\\
&\quad + \mathrm{var}(\mathrm{vec}(Y_0)) + \int_0^h (B\otimes I_d + I_d\otimes B)\,\mathrm{acov}_Y(s)\,\mathrm{d}s\\
&= \mathrm{var}(\mathrm{vec}(Y_0)) + \int_0^h (B\otimes I_d + I_d\otimes B + \sigma_L A\otimes A)\,\mathrm{acov}_Y(s)\,\mathrm{d}s,
\end{aligned}
$$

where we have used a Fubini argument, Lemma 6.9 and $E([L,L]_1^{\mathfrak{d}}) = \sigma_L I_d$. Regarding the use of Lemma 6.9, we observe that $\|(Y_{s-}+C)^{1/2}\otimes(Y_{s-}+C)^{1/2}\|_2 = \|Y_{s-}+C\|_2$ and hence the required local boundedness is ensured by the second-order stationarity of $Y$.

The ordinary differential equation (4.10) is now immediate and to conclude the proof, it suffices to note that $\mathrm{acov}_Y(0) = \mathrm{var}(\mathrm{vec}(Y_0))$ and thus solving the ODE gives

$$
\mathrm{acov}_Y(h) = \mathrm{acov}_V(h) = \mathrm{e}^{(B\otimes I_d + I_d\otimes B + \sigma_L A\otimes A)h}\,\mathrm{var}(\mathrm{vec}(Y_0)), \qquad h\geq 0. \qquad \square
$$

**Proof of Lemma 4.14.** Let $\varepsilon$ be as in the definition of type $\widetilde{G}$ and let $\nu_\varepsilon$ be its Lévy measure. Then, by [2], Proposition 3.1, $L$ has Lévy density $u(x) = \int_{\mathbb{R}^+}\phi_d(x;\tau I_d)\nu_\varepsilon(\mathrm{d}\tau)$, where $\phi_d(\cdot;\Sigma)$ denotes the density of the $d$-dimensional normal distribution with variance

$\Sigma$. Hence,

$$\begin{aligned}E([\operatorname{vec}([L,L]^{\mathfrak{d}}),\operatorname{vec}([L,L]^{\mathfrak{d}})]_1^{\mathfrak{d}}) &= \int_{\mathbb{R}^d}(xx^*)\otimes(xx^*)u(x)\,\mathrm{d}x\\ &= \int_{\mathbb{R}^+}\int_{\mathbb{R}^d}(xx^*)\otimes(xx^*)\phi_d(x;\tau I_d)\,\mathrm{d}x\,\nu_\varepsilon(\mathrm{d}\tau)\\ &= \int_{\mathbb{R}^+}\tau^2\nu_\varepsilon(\mathrm{d}\tau)(I_{d^2}+K_d+\operatorname{vec}(I_d)\operatorname{vec}(I_d)^*),\end{aligned}$$

using [30], Theorem 4.3. Now, set $\rho_L := \int_{\mathbb{R}^+}\tau^2\nu_\varepsilon(\mathrm{d}\tau)$ and note that the finiteness follows from the definition of type $\widetilde{G}$ and the assumed finiteness of the fourth moment of $L$. □

**Proof of Theorem 4.15.** Proposition 4.7 again ensures the existence and local boundedness of the second moment needed in (i).

The definition of quadratic variation (see [5], Lemma 5.11, for a special version in the context of matrix and vector multiplication) implies that

$$\begin{aligned}\operatorname{vec}(Y_t)\operatorname{vec}(Y_t)^* &= \operatorname{vec}(Y_0)\operatorname{vec}(Y_0)^*\\ &\quad+\int_0^t((B\otimes I_d+I_d\otimes B)\operatorname{vec}(Y_{s-})\operatorname{vec}(Y_{s-})^*\\ &\qquad+\operatorname{vec}(Y_{s-})\operatorname{vec}(Y_{s-})^*(B^*\otimes I_d+I_d\otimes B^*))\,\mathrm{d}s\\ &\quad+\int_0^t(A\otimes A)((C+Y_{s-})^{1/2}\otimes(C+Y_{s-})^{1/2})\,\mathrm{d}\operatorname{vec}([L,L]_s^{\mathfrak{d}})\operatorname{vec}(Y_{s-})^*\\ &\quad+\int_0^t\operatorname{vec}(Y_{s-})\,\mathrm{d}\operatorname{vec}([L,L]_s^{\mathfrak{d}})^*((C+Y_{s-})^{1/2}\otimes(C+Y_{s-})^{1/2})(A^*\otimes A^*)\\ &\quad+[\operatorname{vec}(Y),\operatorname{vec}(Y)]_t.\end{aligned}$$

Moreover, setting $\mathcal{V}_t = ((C+Y_t)^{1/2}\otimes(C+Y_t)^{1/2})$ we obtain

$$\begin{aligned}[\operatorname{vec}(Y),\operatorname{vec}(Y)^*]_t &= \int_0^t(A\otimes A)\mathcal{V}_{s-}\,\mathrm{d}([\operatorname{vec}([L,L]^{\mathfrak{d}}),\operatorname{vec}([L,L]^{\mathfrak{d}})]_s^{\mathfrak{d}})\mathcal{V}_{s-}(A^*\otimes A^*)\\ &= \int_0^t\int_{\mathbb{R}^d}(A\otimes A)\mathcal{V}_{s-}\operatorname{vec}(xx^*)\operatorname{vec}(xx^*)^*\mathcal{V}_{s-}(A^*\otimes A^*)\mu_L(\mathrm{d}s,\mathrm{d}x),\end{aligned}$$

since $\operatorname{vec}(Y_t)$ is the sum of an absolutely continuous component and a pure jump process of finite variation.

Using a Fubini argument, the stochastic continuity, Lemma 6.9 and the Assumptions 4.2, 4.3 made on the moments of $\nu_L$, we obtain

$$E(\operatorname{vec}(Y_t)\operatorname{vec}(Y_t)^*)$$

$$
\begin{aligned}
&= E(\operatorname{vec}(Y_0)\operatorname{vec}(Y_0)^*)\\
&\quad + \int_0^t ((B\otimes I_d + I_d\otimes B)E(\operatorname{vec}(Y_s)\operatorname{vec}(Y_s)^*)\\
&\qquad\quad + E(\operatorname{vec}(Y_s)\operatorname{vec}(Y_s)^*)(B^*\otimes I_d + I_d\otimes B^*))\,\mathrm{d}s\\
&\quad + \sigma_L\int_0^t (A\otimes A)E(\operatorname{vec}(C+Y_s)\operatorname{vec}(Y_s)^*)\,\mathrm{d}s\\
&\quad + \sigma_L\int_0^t E(\operatorname{vec}(Y_s)\operatorname{vec}(C+Y_s)^*)(A^*\otimes A^*)\,\mathrm{d}s\\
&\quad + \int_0^t (A\otimes A)E(\mathcal{V}_{s-}\rho_L(I_{d^2}+K_d+\operatorname{vec}(I_d)\operatorname{vec}(I_d)^*)\mathcal{V}_{s-})(A^*\otimes A^*)\,\mathrm{d}s.
\end{aligned}
$$

With the definition of $\mathcal{V}_t$, it follows that

$$
\begin{aligned}
E(\mathcal{V}_{s-}I_{d^2}\mathcal{V}_{s-}) &= E(\mathcal{V}_s^2) = E((C+Y_s)\otimes(C+Y_{s-}))\\
&= \mathbf{Q}E(\operatorname{vec}(C+Y_s)\operatorname{vec}(C+Y_s)^*),\\
E(\mathcal{V}_{s-}\operatorname{vec}(I_d)\operatorname{vec}(I_d)^*\mathcal{V}_{s-}) &= E(\operatorname{vec}(C+Y_s)\operatorname{vec}(C+Y_s)^*),\\
E(\mathcal{V}_{s-}K_d\mathcal{V}_{s-}) &= K_d E((C+Y_s)\otimes(C+Y_s))\\
&= K_d\mathbf{Q}E(\operatorname{vec}(C+Y_s)\operatorname{vec}(C+Y_s)^*),
\end{aligned}
$$

using [30], Theorem 3.1 (xii), in the last identity. Inserting these formulae into the above result and noting that in the stationary case, the integrands need to sum to zero gives (4.17). Vectorizing then immediately establishes (4.18).

Likewise, we obtain (4.16) in the non-stationary case by inserting the formulae above, vectorizing and differentiating. □

**Proof of Lemma 4.17.** We have

$$
\begin{aligned}
&\rho_L\|\operatorname{vec}(I_{d^2}+K_d+\operatorname{vec}(I_d)\operatorname{vec}(I_d)^*)\|_{\widetilde{B,S}}\\
&\quad = \left\|\int_{\mathbb{R}^d}\operatorname{vec}((xx^*)\otimes(xx^*))\nu_L(\mathrm{d}x)\right\|_{\widetilde{B,S}} \qquad (6.7)\\
&\quad \le \int_{\mathbb{R}^d}\|\operatorname{vec}((xx^*)\otimes(xx^*))\|_{\widetilde{B,S}}\nu_L(\mathrm{d}x) = \int_{\mathbb{R}^d}\|\operatorname{vec}(xx^*)\|_{B,S}^2\nu_L(\mathrm{d}x)
\end{aligned}
$$

since the definition of $\|\cdot\|_{\widetilde{B,S}}$ implies that

$$
\begin{aligned}
\|\operatorname{vec}((xx^*)\otimes(xx^*))\|_{\widetilde{B,S}} &= \|\mathcal{S}^{-1}(x\otimes x)\otimes(x\otimes x)\|_2\\
&= \|(S^{-1}\otimes S^{-1})(x\otimes x)\|_2^2 = \|\operatorname{vec}(xx^*)\|_{B,S}^2,
\end{aligned}
$$

using the fact that $\|z \otimes z\|_2 = \|z\|_2^2$ for all $z \in \mathbb{R}^{d^2}$.

For $k = 2$, the condition (4.4) becomes

$$2\|S\|_2^2\|S^{-1}\|_2^2 K_{2,B}\|A \otimes A\|_{B,S} \int_{\mathbb{R}^d} \|\operatorname{vec}(xx^*)\|_{B,S}\nu_L(\mathrm{d}x)$$
$$+\|S\|_2^4\|S^{-1}\|_2^4 K_{2,B}^2\|A \otimes A\|_{B,S}^2 \int_{\mathbb{R}^d} \|\operatorname{vec}(xx^*)\|_{B,S}^2\nu_L(\mathrm{d}x) < -4\lambda.$$

Using (6.7) and results from the proof of Lemma 4.10 gives

$$2\sigma_L\|A \otimes A\|_{B,S} + K_{2,B}^2\|A \otimes A\|_{B,S}^2\rho_L\|\operatorname{vec}(I_{d^2} + K_d + \operatorname{vec}(I_d)\operatorname{vec}(I_d)^*)\|_{\widetilde{B,S}} < -4\lambda.$$

Combining $\|(A \otimes A) \otimes I_{d^2} + I_{d^2} \otimes (A \otimes A)\|_{\widetilde{B,S}} \leq 2\|A \otimes A\|_{B,S}$ and $\|A \otimes A \otimes A \otimes A\|_{\widetilde{B,S}} = \|A \otimes A\|_{B,S}^2$, which are elementary to prove, with (4.19) leads to $\|\mathscr{C} - \mathscr{B} \otimes I_{d^2} - I_{d^2} \otimes \mathscr{B}\|_{\widetilde{B,S}} < -4\lambda$.

Since $\mathcal{S}^{-1}(\mathscr{B} \otimes I_{d^2} + I_{d^2} \otimes \mathscr{B})\mathcal{S}$ is diagonal and $\max(\Re(\sigma(\mathscr{B} \otimes I_{d^2} + I_{d^2} \otimes \mathscr{B}))) = 4\lambda$, the Bauer–Fike theorem (see [23], Theorem 6.3.2 and its proof) and arguments as in the proof of Lemma 4.10 complete this proof. □

**Proof of Lemma 4.18.** That $S$ is unitary implies that $K_{2,B} = 1$ and all the norms used are actually the Euclidean norm or the operator norm induced by it. Hence, we have to show that

$$\|\mathcal{Q} + \mathcal{K}_d\mathcal{Q} + I_{d^4}\|_2 \leq \|\operatorname{vec}(I_{d^2} + K_d + \operatorname{vec}(I_d)\operatorname{vec}(I_d)^*)\|_2.$$

For $d = 1$, one calculates both sides to be equal to 3.

In general, we know from the fact that $\mathcal{K}_d$ and $\mathcal{Q}$ are permutation matrices that the operator norms are 1. Hence, $\|\mathcal{Q} + \mathcal{K}_d\mathcal{Q} + I_{d^4}\|_2 \leq 3$. Furthermore, the entries of $K_d$ and $\operatorname{vec}(I_d)\operatorname{vec}(I_d)^*$ are either 1 or 0. Therefore, $\|\operatorname{vec}(I_{d^2} + K_d + \operatorname{vec}(I_d)\operatorname{vec}(I_d)^*)\|_2 \geq \|\operatorname{vec}(I_{d^2})\|_2 = d$. This shows the inequality for $d \geq 3$.

In the remaining case $d = 2$, we have

$$\|\operatorname{vec}(I_{d^2} + K_d + \operatorname{vec}(I_d)\operatorname{vec}(I_d)^*)\|_2 \geq \|\operatorname{vec}(I_{d^2} + \operatorname{vec}(I_d)\operatorname{vec}(I_d)^*)\|_2 = \sqrt{12} > 3,$$

which again establishes the claimed inequality. □

To prove the asymptotic second-order stationarity, we need the following general lemma on differential equations which is elementary, but not to be found in the literature, to the best of our knowledge.

**Lemma 6.10.** *Let $f : \mathbb{R}^+ \to \mathbb{R}^d$ be continuous and $A \in M_d(\mathbb{R})$ with $\sigma(A) \subset (-\infty, 0) + \mathrm{i}\mathbb{R}$.*

*If $\lim_{t\to\infty} f(t) = \xi$ with $\xi \in \mathbb{R}^d$, then, for any initial value $x_0 \in \mathbb{R}^d$, the solution $x$ to the differential equation*

$$\frac{\mathrm{d}x(t)}{\mathrm{d}t} = Ax(t) + f(t)$$

*satisfies $\lim_{t\to\infty} x(t) = -A^{-1}\xi$.*

**Proof.** It holds that

$$x(t) = \mathrm{e}^{At}x_0 + \int_0^t \mathrm{e}^{A(t-s)} f(s)\,\mathrm{d}s.$$

Since

$$\lim_{t\to\infty} \int_0^t \mathrm{e}^{A(t-s)}\xi\,\mathrm{d}s = -A^{-1}\xi \quad \text{and} \quad \lim_{t\to\infty} \mathrm{e}^{At}x_0 = 0,$$

it suffices to show that $\lim_{t\to\infty} \|\int_0^t \mathrm{e}^{A(t-s)}(f(s)-\xi)\,\mathrm{d}s\| = 0$. Fix $\varepsilon > 0$. There exist $t^*, t^{**} > 0$ with $t^* \le t^{**}$ such that $\|f(t)-\xi\| < \varepsilon$ for all $t \ge t^*$ and

$$\left\| \int_0^{t^*} \mathrm{e}^{A(t-s)}(f(s)-\xi)\,\mathrm{d}s \right\| < \varepsilon \qquad \text{for all } t \ge t^{**}.$$

Hence,

$$\left\| \int_0^t \mathrm{e}^{A(t-s)}(f(s)-\xi)\,\mathrm{d}s \right\| \le \left(1 + \int_{t^*}^t \|\mathrm{e}^{A(t-s)}\|\,\mathrm{d}s\right)\varepsilon \le \left(1 + \int_0^\infty \|\mathrm{e}^{As}\|\,\mathrm{d}s\right)\varepsilon \qquad \forall t \ge t^{**}.$$

Since the last integral is finite and $\varepsilon$ was arbitrary, this completes the proof. □

**Proof of Theorem 4.20.** (i) follows from Proposition 4.7 and Theorems 4.8, 4.11 and 4.15.

Regarding (ii), Proposition 4.7 ensures that $E(\|Y_t\|^2) < \infty$ for all $t \in \mathbb{R}^+$. The convergence of the expectation has already been noted in Remark 4.9 and the convergence of the variance follows from (4.16) and the previous lemma. (4.12) then implies the convergence of the autocovariance. □

## 6.3. Proofs for Section 5

**Proof of Proposition 5.2.** Proposition 4.7 ensures that $E(\|Y_t\|)$ and hence $E(\|V_t\|)$ is finite and locally bounded. Since $E(\|V_t^{1/2}\|_2^2) = E(\|V_t\|_2)$, the standard $L^2$-stochastic integration theory (see [1], Section 4.2.1, for example) establishes that $E(\|G_t\|^2) < \infty$ for all $t \in \mathbb{R}^+$.

If we now let Assumptions 5.1 and 5.2 be satisfied, then (5.2) is clear and (5.5) is a straightforward consequence of the Itô isometry. The latter also implies that

$$\begin{aligned}\operatorname{var}(\mathbf{G}_1) &= E(\mathbf{G}_1\mathbf{G}_1^*) = E\left(\int_0^\Delta V_{s-}^{1/2} E(L_1L_1^*) V_{s-}^{1/2}\,\mathrm{d}s\right) = (\sigma_L + \sigma_W) E\left(\int_0^\Delta V_{s-}\,\mathrm{d}s\right)\\ &= (\sigma_L + \sigma_W)\Delta E(V_0).\end{aligned}$$

□

**Proof of Proposition 5.4.** (i) We first show the finiteness of the moments. Using the Euclidean norm and its operator norm, it is clear that $E(\|G_t\|^4) < \infty$ if and only if $E(\|G_t G_t^*\|^2) < \infty$. Denoting the $d$ components of $G$ by $G_i$ with $i = 1, \dots, d$, we have that $E(\|G_t\|^4) < \infty$ is equivalent to $E(|G_{i,t}|^4) < \infty$ for all $i \in \{1, \dots, d\}$. But the Burkholder–Davis–Gundy inequalities (see [36], page 222) give that $E(|G_{i,t}|^4) < \infty$ provided $E([G_i, G_i]_t^2) < \infty$. The latter is, in turn, ensured by $E(\|[G, G]_t\|^2) < \infty$ simultaneously for all $i \in \{1, \dots, d\}$.

Next, we observe that $[L, L]_t = \tau_L t + \int_0^t xx^* \mu_L(\mathrm{d}s, \mathrm{d}x)$ with $\tau_L \in \mathbb{S}_d^+$ and that, moreover, $[G, G]_t = \int_0^t V_{s-}^{1/2} \,\mathrm{d}[L, L]_s \, V_{s-}^{1/2}$ or, equivalently, $\mathrm{vec}([G, G]_t) = \int_0^t V_{s-}^{1/2} \otimes V_{s-}^{1/2} \,\mathrm{d}\,\mathrm{vec}([L, L]_s)$. Since Proposition 4.7 ensures the finiteness and local boundedness of $E(\|V_s\|_2^2) = E(\|V_s \otimes V_s\|_2)$ and $E(\|L_1\|^4)$ the finiteness of $E(\|[L, L]_1\|^2)$, the standard $L^2$-stochastic integration theory immediately gives $E(\|[G, G]_t\|^2) < \infty$.

(ii) If we now let Assumptions 5.1 and 5.2 be satisfied, then (5.6) has already been shown in the last proposition.

It thus remains to establish (5.7). Let $h \in \mathbb{N}$. The definition of the quadratic (co)variation (see [5], Lemma 5.11, in particular) implies that

$$\begin{aligned}\mathbf{G}_{h+1}\mathbf{G}_{h+1}^* = &\int_{h\Delta}^{(h+1)\Delta} V_{s-}^{1/2} \,\mathrm{d}L_s \left(\int_{h\Delta}^{s} \mathrm{d}L_u^* \, V_{u-}^{1/2}\right) + \int_{h\Delta}^{(h+1)\Delta} \left(\int_{h\Delta}^{s} V_{u-}^{1/2} \,\mathrm{d}L_u\right) \mathrm{d}L_s^* \, V_{s-}^{1/2} \\ &+ \left[\int_{h\Delta}^{(h+1)\Delta} V_{s-}^{1/2} \,\mathrm{d}L_s, \int_{h\Delta}^{(h+1)\Delta} V_{s-}^{1/2} \,\mathrm{d}L_s\right]\end{aligned}$$

with

$$\begin{aligned}&\left[\int_{h\Delta}^{(h+1)\Delta} V_{s-}^{1/2} \,\mathrm{d}L_s, \int_{h\Delta}^{(h+1)\Delta} V_{s-}^{1/2} \,\mathrm{d}L_s\right] \\ &\quad = \int_{h\Delta}^{(h+1)\Delta} V_{s-}^{1/2} \,\mathrm{d}[L, L]_s \, V_{s-}^{1/2} \\ &\quad = \sigma_W \int_{h\Delta}^{(h+1)\Delta} V_{s-} \,\mathrm{d}s + \int_{h\Delta}^{(h+1)\Delta} V_{s-}^{1/2} \,\mathrm{d}[L, L]_s^{\mathfrak{d}} \, V_{s-}^{1/2}.\end{aligned}$$

We now condition upon $\mathcal{F}_\Delta$ and denote by $Y(y, (L_r - L_{t_0})_{r \geq t_0}, t_0, t)$ with $t_0 \in \mathbb{R}^+$ the solution of

$$\mathrm{d}Y_t = (BY_{t-} + Y_{t-}B^*) \,\mathrm{d}t + A(Y_{t-} + C)^{1/2} \,\mathrm{d}[L, L]_t^{\mathfrak{d}} \, (Y_{t-} + C)^{1/2} A^*$$

for $t \geq t_0$ with $Y_{t_0} = y$. Furthermore, we denote by $E_{L,t_0}(\cdot)$ the expectation taken with respect to $(L_r - L_{t_0})_{r \geq t_0}$ only. Using Theorem 4.4, $E(L_1) = 0$ and the fact that the increments of $(L_r)_{r \geq \Delta}$ are independent of $\mathcal{F}_\Delta$, one obtains

$$\begin{aligned}&E\left(\int_{h\Delta}^{(h+1)\Delta} V_{s-}^{1/2} \,\mathrm{d}L_s \left(\int_{h\Delta}^{s} \mathrm{d}L_u^* V_{u-}^{1/2}\right) \Big| \mathcal{F}_\Delta\right) \\ &\quad = E\left(\int_{h\Delta}^{(h+1)\Delta} V_{s-}^{1/2} \,\mathrm{d}L_s \left(\int_{h\Delta}^{s} \mathrm{d}L_u^* V_{u-}^{1/2}\right) \Big| Y_\Delta\right)\end{aligned}$$

$$= E_{L,\Delta}\left(\int_{h\Delta}^{(h+1)\Delta}\int_{h\Delta}^{s} V_{Y_\Delta,s-}^{1/2}\,\mathrm{d}L_s\,\mathrm{d}L_u^*\,V_{Y_\Delta,u-}^{1/2}\right) = 0,$$

where $V_{Y_\Delta,t} := Y(Y_\Delta,(L_r - L_\Delta)_{r\geq\Delta},\Delta,t) + C$ and, likewise,

$$E\left(\int_{h\Delta}^{(h+1)\Delta}\left(\int_{h\Delta}^{s} V_{u-}^{1/2}\,\mathrm{d}L_u\right)\mathrm{d}L_s^* V_{s-}^{1/2}\Big|\mathcal{F}_\Delta\right) = 0.$$

Moreover, using the moment assumptions and the compensation formula, we have

$$E\left(\sigma_W \int_{h\Delta}^{(h+1)\Delta} V_{s-}\,\mathrm{d}s + \int_{h\Delta}^{(h+1)\Delta} V_{s-}^{1/2}\,\mathrm{d}[L,L]_s^{\mathfrak{d}} V_{s-}^{1/2}\Big|\mathcal{F}_\Delta\right)$$
$$= (\sigma_L + \sigma_W)\int_{h\Delta}^{(h+1)\Delta} E_{L,\Delta}(V_{Y_\Delta,s})\,\mathrm{d}s.$$

Equation (4.7) implies that

$$\begin{aligned}
&\int_{h\Delta}^{(h+1)\Delta} E_{L,\Delta}(\mathrm{vec}(V_{Y_\Delta,s}))\,\mathrm{d}s\\
&= \int_{h\Delta}^{(h+1)\Delta} (\mathrm{vec}(C) + \mathrm{e}^{\mathscr{B}(s-\Delta)}(\mathrm{vec}(Y_\Delta) + \sigma_L \mathscr{B}^{-1}(A\otimes A)\,\mathrm{vec}(C)))\,\mathrm{d}s\\
&\quad - \int_{h\Delta}^{(h+1)\Delta} \sigma_L \mathscr{B}^{-1}(A\otimes A)\,\mathrm{vec}(C)\,\mathrm{d}s\\
&= \Delta\,\mathrm{vec}(C) + \mathscr{B}^{-1}\mathrm{e}^{\mathscr{B}\Delta h}(I_{d^2} - \mathrm{e}^{-\mathscr{B}\Delta})(\mathrm{vec}(Y_\Delta) - E(\mathrm{vec}(Y_0))) + \Delta E(\mathrm{vec}(Y_0)).
\end{aligned}$$

Combining the above results, we get

$$\begin{aligned}
&E(\mathrm{vec}(\mathbf{G}_{h+1}\mathbf{G}_{h+1}^*)(\mathrm{vec}(\mathbf{G}_1\mathbf{G}_1^*))^*)\\
&\quad = E(E(\mathrm{vec}(\mathbf{G}_{h+1}\mathbf{G}_{h+1}^*)|\mathcal{F}_\Delta)(\mathrm{vec}(\mathbf{G}_1\mathbf{G}_1^*))^*)\\
&\quad = E(\mathrm{vec}(\mathbf{G}_1\mathbf{G}_1^*))(E(\mathrm{vec}(\mathbf{G}_1\mathbf{G}_1^*)))^* + (\sigma_L+\sigma_W)\mathrm{e}^{\mathscr{B}\Delta h}\mathscr{B}^{-1}(I_{d^2} - \mathrm{e}^{-\mathscr{B}\Delta})\\
&\qquad \times (E(\mathrm{vec}(Y_\Delta)(\mathrm{vec}(\mathbf{G}_1\mathbf{G}_1^*))^*) - E(\mathrm{vec}(Y_0))(E(\mathrm{vec}(\mathbf{G}_1\mathbf{G}_1^*)))^*).
\end{aligned}$$

Using the stationarity of $Y$, this establishes (5.7). □

## Acknowledgements

The author is very grateful to Jean Jacod and his Ph.D. advisor Claudia Klüppelberg for very helpful comments and suggestions. Furthermore, financial support by the Deutsche Forschungsgemeinschaft through the graduate programme "Angewandte Algorithmische Mathematik" at the Munich University of Technology is gratefully acknowledged.